\documentclass[12pt]{article}
\usepackage{theorem}
\usepackage{amssymb}
\usepackage{graphicx}
\usepackage[mathscr]{eucal}
\newtheorem{prop}{}[section]

{\theorembodyfont{\upshape} \newtheorem{rema}[prop]{}}
\newcommand{\boma}[1]{{\mbox{\boldmath $#1$} }}

\begin{document}
\newcommand{\binom}[2]{\left( \barray{c} #1 \\ #2 \farray \right)}
\newcommand{\uper}[1]{\stackrel{\barray{c} {~} \\ \mbox{\footnotesize{#1}}\farray}{\longrightarrow} }
\newcommand{\nop}[1]{ \|#1\|_{\piu} }
\newcommand{\no}[1]{ \|#1\| }
\newcommand{\nom}[1]{ \|#1\|_{\meno} }
\newcommand{\UU}[1]{e^{#1 \AA}}
\newcommand{\UD}[1]{e^{#1 \Delta}}
\newcommand{\bb}[1]{\mathbb{{#1}}}
\newcommand{\HO}[1]{\bb{H}^{{#1}}}
\newcommand{\Hz}[1]{\bb{H}^{{#1}}_{\zz}}
\newcommand{\Hs}[1]{\bb{H}^{{#1}}_{\ss}}
\newcommand{\Hg}[1]{\bb{H}^{{#1}}_{\gam}}
\newcommand{\HM}[1]{\bb{H}^{{#1}}_{\so}}
\newcommand{\vers}[1]{\widehat{#1}}
\def\tvainf{\vspace{-0.4cm} \barray{ccc} \vspace{-0,1cm}{~}
\\ \vspace{-0.2cm} \longrightarrow \\ \vspace{-0.2cm} \scriptstyle{T \vain + \infty} \farray}
\def\ua{u_{\tt{a}}}
\def\Ta{T_{\tt{a}}}
\def\Tc{T_{\tt{c}}}
\def\tu{\xi}
\def\ep{\xi}
\def\iep{(-1,+\infty)}
\def\Do{\mathscr{E}}
\def\CA{\mathcal{C}}
\def\CB{\mathcal{D}}
\def\ca{c}
\def\cb{d}
\def\op{\,\mbox{\scriptsize{or}}\,}
\def\er{\epsilon}
\def\erd{\er_0}
\def\vk{\vers{k}}
\def\vh{\vers{h}}
\def\Qn{\mathfrak{K}_n}
\def\Ed{\hat{E}}
\def\um{u_{-}}
\def\up{u_{+}}
\def\el{t}
\def\em{z}
\def\uu{\lambda}
\def\dK{\delta {\mathscr K}}
\def\dG{\delta {\mathscr G}}
\def\DK{\Delta {\mathscr K}}
\def\DG{\Delta {\mathscr G}}
\def\Ti{{\mathscr T}}
\def\Km{{\mathscr K}}
\def\Ll{\mathscr{L}}
\def\Hh{\mathscr{H}}
\def\Mm{{\mathscr M}}
\def\Nn{{\mathscr N}}
\def\Rr{{\mathscr R}}
\def\Gg{{\mathscr G}}
\def\Zz{Z}
\def\Ss{{\mathscr S}}
\def\Fe{{\mathscr F}}
\def\Ei{{\mathscr E}}
\def\Ww{{\mathscr W}}
\def\we{\wedge}
\def\We{\bigwedge}
\def\dbar{\hat{d}}
\def\Cc{\mathscr{C}}
\def\SZ{\mathcal{S}}
\def\TZ{\mathfrak{S}}
\def\CQ{S}
\def\GQ{T}
\def\C{{C_d}}
\def\Ac{\overline{A}}
\def\Bc{\overline{B}}
\def\Xd{\Zd_{0 k} \cap \Bd_{0 k}(\ro)}
\def\Yd{\Zd_{0 k} \setminus \Bd_{0 k}(\ro) }
\def\comple{\scriptscriptstyle{\complessi}}
\def\nume{0.407}
\def\numerob{0.00724}
\def\deln{7/10}
\def\delnn{\dd{7 \over 10}}
\def\e{c}
\def\p{p}
\def\z{z}
\def\symd{{\mathfrak S}_d}
\def\del{\omega}
\def\Del{\delta}
\def\Di{\Delta}
\def\mmu{\hat{\mu}}
\def\rot{\mbox{rot}\,}
\def\curl{\mbox{curl}\,}
\def\XS{\boma{x}}
\def\TS{\boma{t}}
\def\Lam{\boma{\eta}}
\def\DS{\boma{\rho}}
\def\KS{\boma{k}}
\def\LS{\boma{\lambda}}
\def\PR{\boma{p}}
\def\VS{\boma{v}}
\def\ski{\! \! \! \! \! \! \! \! \! \! \! \! \! \!}
\def\h{L}
\def\EM{M}
\def\EMP{M'}
\def\R{R}
\def\E{E}
\def\FFf{\mathscr{F}}
\def\A{F}
\def\Xim{\Xi_{\meno}}
\def\Ximn{\Xi_{n-1}}
\def\lan{\lambda}
\def\om{\omega}
\def\Om{\Omega}
\def\Sim{\Sigm}
\def\Sip{\Delta \Sigm}
\def\Sigm{{\mathscr{S}}}
\def\Ki{{\mathscr{K}}}
\def\Hi{{\mathscr{H}}}
\def\zz{{\scriptscriptstyle{0}}}
\def\ss{{\scriptscriptstyle{\Sigma}}}
\def\gam{{\scriptscriptstyle{\Gamma}}}
\def\so{\ss \zz}
\def\Dz{\bb{\DD}'_{\zz}}
\def\Ds{\bb{\DD}'_{\ss}}
\def\Dsz{\bb{\DD}'_{\so}}
\def\Dg{\bb{\DD}'_{\gam}}
\def\Ls{\bb{L}^2_{\ss}}
\def\Lg{\bb{L}^2_{\gam}}
\def\bF{{\bb{V}}}
\def\Fz{\bF_{\zz}}
\def\Fs{\bF_\ss}
\def\Fg{\bF_\gam}
\def\Pre{P}
\def\UUU{{\mathcal U}}
\def\fiapp{\phi}
\def\PU{P1}
\def\PD{P2}
\def\PT{P3}
\def\PQ{P4}
\def\PC{P5}
\def\PS{P6}
\def\Q{P6}
\def\X{Q2}
\def\Xp{Q3}
\def\Vi{V}
\def\bVi{\bb{V}}
\def\K{V}
\def\Ks{\bb{\K}_\ss}
\def\Kz{\bb{\K}_0}
\def\KM{\bb{\K}_{\, \so}}
\def\HGG{\bb{H}^\G}
\def\HG{\bb{H}^\G_{\so}}
\def\EG{{\mathfrak{P}}^{\G}}
\def\G{G}
\def\de{\delta}
\def\esp{\sigma}
\def\dd{\displaystyle}
\def\LP{\mathfrak{L}}
\def\dive{\mbox{div}\,}
\def\la{\langle}
\def\ra{\rangle}
\def\um{u_{\meno}}
\def\uv{\mu_{\meno}}
\def\Fp{ {\textbf F_{\piu}} }
\def\Ff{ {\textbf F} }
\def\Fm{ {\textbf F_{\meno}} }
\def\piu{\scriptscriptstyle{+}}
\def\meno{\scriptscriptstyle{-}}
\def\omeno{\scriptscriptstyle{\ominus}}
\def\Tt{ {\mathscr T} }
\def\Xx{ {\textbf X} }
\def\Yy{ {\textbf Y} }
\def\VP{{\mbox{\tt VP}}}
\def\CP{{\mbox{\tt CP}}}
\def\cp{$\CP(f_0, t_0)\,$}
\def\cop{$\CP(f_0)\,$}
\def\copn{$\CP_n(f_0)\,$}
\def\vp{$\VP(f_0, t_0)\,$}
\def\vop{$\VP(f_0)\,$}
\def\vopn{$\VP_n(f_0)\,$}
\def\vopdue{$\VP_2(f_0)\,$}
\def\leqs{\leqslant}
\def\geqs{\geqslant}
\def\mat{{\frak g}}
\def\tG{t_{\scriptscriptstyle{G}}}
\def\tN{t_{\scriptscriptstyle{N}}}
\def\TK{t_{\scriptscriptstyle{K}}}
\def\CK{C_{\scriptscriptstyle{K}}}
\def\CN{C_{\scriptscriptstyle{N}}}
\def\CG{C_{\scriptscriptstyle{G}}}
\def\CCG{{\mathscr{C}}_{\scriptscriptstyle{G}}}
\def\tf{{\tt f}}
\def\ti{{\tt t}}
\def\ta{{\tt a}}
\def\tc{{\tt c}}
\def\tF{{\tt R}}
\def\P{{\mathscr P}}
\def\V{{\mathscr V}}
\def\TI{\tilde{I}}
\def\TJ{\tilde{J}}
\def\Lin{\mbox{Lin}}
\def\Hinfc{ H^{\infty}(\reali^d, \complessi) }
\def\Hnc{ H^{n}(\reali^d, \complessi) }
\def\Hmc{ H^{m}(\reali^d, \complessi) }
\def\Hac{ H^{a}(\reali^d, \complessi) }
\def\Dc{\DD(\reali^d, \complessi)}
\def\Dpc{\DD'(\reali^d, \complessi)}
\def\Sc{\SS(\reali^d, \complessi)}
\def\Spc{\SS'(\reali^d, \complessi)}
\def\Ldc{L^{2}(\reali^d, \complessi)}
\def\Lpc{L^{p}(\reali^d, \complessi)}
\def\Lqc{L^{q}(\reali^d, \complessi)}
\def\Lrc{L^{r}(\reali^d, \complessi)}
\def\Hinfr{ H^{\infty}(\reali^d, \reali) }
\def\Hnr{ H^{n}(\reali^d, \reali) }
\def\Hmr{ H^{m}(\reali^d, \reali) }
\def\Har{ H^{a}(\reali^d, \reali) }
\def\Dr{\DD(\reali^d, \reali)}
\def\Dpr{\DD'(\reali^d, \reali)}
\def\Sr{\SS(\reali^d, \reali)}
\def\Spr{\SS'(\reali^d, \reali)}
\def\Ldr{L^{2}(\reali^d, \reali)}
\def\Hinfk{ H^{\infty}(\reali^d, \KKK) }
\def\Hnk{ H^{n}(\reali^d, \KKK) }
\def\Hmk{ H^{m}(\reali^d, \KKK) }
\def\Hak{ H^{a}(\reali^d, \KKK) }
\def\Dk{\DD(\reali^d, \KKK)}
\def\Dpk{\DD'(\reali^d, \KKK)}
\def\Sk{\SS(\reali^d, \KKK)}
\def\Spk{\SS'(\reali^d, \KKK)}
\def\Ldk{L^{2}(\reali^d, \KKK)}
\def\Knb{K^{best}_n}
\def\sc{\cdot}
\def\k{\mbox{{\tt k}}}
\def\x{\mbox{{\tt x}}}
\def\g{ {\textbf g} }
\def\QQQ{ {\textbf Q} }
\def\AAA{ {\textbf A} }
\def\gr{\mbox{gr}}
\def\sgr{\mbox{sgr}}
\def\loc{\mbox{loc}}
\def\PZ{{\Lambda}}
\def\PZAL{\mbox{P}^{0}_\alpha}
\def\epsilona{\epsilon^{\scriptscriptstyle{<}}}
\def\epsilonb{\epsilon^{\scriptscriptstyle{>}}}
\def\lgraffa{ \mbox{\Large $\{$ } \hskip -0.2cm}
\def\rgraffa{ \mbox{\Large $\}$ } }
\def\restriction{\upharpoonright}
\def\M{{\scriptscriptstyle{M}}}
\def\m{m}
\def\Fre{Fr\'echet~}
\def\I{{\mathcal N}}
\def\ap{{\scriptscriptstyle{ap}}}
\def\fiap{\varphi_{\ap}}
\def\dfiap{{\dot \varphi}_{\ap}}
\def\DDD{ {\mathfrak D} }
\def\BBB{ {\textbf B} }
\def\EEE{ {\textbf E} }
\def\GGG{ {\textbf G} }
\def\TTT{ {\textbf T} }
\def\KKK{ {\textbf K} }
\def\HHH{ {\textbf K} }
\def\FFi{ {\bf \Phi} }
\def\GGam{ {\bf \Gamma} }
\def\sc{ {\scriptstyle{\bullet} }}
\def\a{a}
\def\c{\kappa}
\def\parn{\par\noindent}
\def\teta{M}
\def\elle{L}
\def\ro{\rho}
\def\al{\alpha}
\def\alc{\overline{\al}}
\def\dal{\mathfrak{a}}
\def\si{\sigma}
\def\be{\beta}
\def\dbe{\mathfrak{b}}
\def\bec{\overline{\be}}
\def\dbec{\overline{\dbe}}
\def\ga{\gamma}
\def\tet{\vartheta}
\def\teta{\theta}
\def\ch{\chi}
\def\et{\eta}
\def\complessi{{\bf C}}
\def\len{{\bf L}}
\def\reali{{\bf R}}
\def\interi{{\bf Z}}
\def\Z{{\bf Z}}
\def\naturali{{\bf N}}
\def\To{ {\bf T} }
\def\Td{ {\To}^d }
\def\Tt{ {\To}^3 }
\def\Bd{B^d}
\def\Zd{ \interi^d }
\def\Zt{ \interi^3 }
\def\Zet{{\mathscr{Z}}}
\def\Ze{\Zet^d}
\def\T1{{\textbf To}^{1}}
\def\Sfe{ {\bf S} }
\def\Sd{\Sfe^{d-1}}
\def\St{\Sfe^{2}}
\def\es{s}
\def\FF{\mathcal F}
\def\FFu{ {\textbf F_{1}} }
\def\FFd{ {\textbf F_{2}} }
\def\GG{{\mathcal G} }
\def\EE{{\mathcal E}}
\def\KK{{\mathcal K}}
\def\PP{{\mathcal P}}
\def\PPP{{\mathscr P}}
\def\PN{{\mathcal P}}
\def\PPN{{\mathscr P}}
\def\QQ{{\mathcal Q}}
\def\J{J}
\def\Np{{\hat{N}}}
\def\Lp{{\hat{L}}}
\def\Jp{{\hat{J}}}
\def\Vp{{\hat{V}}}
\def\Ep{{\hat{E}}}
\def\Gp{{\hat{G}}}
\def\Kp{{\hat{K}}}
\def\Ip{{\hat{I}}}
\def\Tp{{\hat{T}}}
\def\Mp{{\hat{M}}}
\def\La{\Lambda}
\def\Ga{\Gamma}
\def\Si{\Sigma}
\def\Upsi{\Upsilon}
\def\Gam{\Gamma}
\def\Gag{{\check{\Gamma}}}
\def\Lap{{\hat{\Lambda}}}
\def\Upsig{{\check{\Upsilon}}}
\def\Kg{{\check{K}}}
\def\ellp{{\hat{\ell}}}
\def\j{j}
\def\jp{{\hat{j}}}
\def\BB{{\mathcal B}}
\def\LL{{\mathcal L}}
\def\MM{{\mathcal U}}
\def\SS{{\mathcal S}}
\def\DD{D}
\def\Dd{{\mathcal D}}
\def\VV{{\mathcal V}}
\def\WW{{\mathcal W}}
\def\OO{{\mathcal O}}
\def\RR{{\mathcal R}}
\def\TT{{\mathcal T}}
\def\AA{{\mathcal A}}
\def\CC{{\mathcal C}}
\def\JJ{{\mathcal J}}
\def\NN{{\mathcal N}}
\def\HH{{\mathcal H}}
\def\XX{{\mathcal X}}
\def\XXX{{\mathscr X}}
\def\YY{{\mathcal Y}}
\def\ZZ{{\mathcal Z}}
\def\cir{{\scriptscriptstyle \circ}}
\def\circa{\thickapprox}
\def\vain{\rightarrow}
\def\parn{\par \noindent}
\def\salto{\vskip 0.2truecm \noindent}
\def\spazio{\vskip 0.5truecm \noindent}
\def\vs1{\vskip 1cm \noindent}
\def\fine{\hfill $\square$ \vskip 0.2cm \noindent}
\def\ffine{\hfill $\lozenge$ \vskip 0.2cm \noindent}
\newcommand{\rref}[1]{(\ref{#1})}
\def\beq{\begin{equation}}
\def\feq{\end{equation}}
\def\beqq{\begin{eqnarray}}
\def\feqq{\end{eqnarray}}
\def\barray{\begin{array}}
\def\farray{\end{array}}
\makeatletter \@addtoreset{equation}{section}
\renewcommand{\theequation}{\thesection.\arabic{equation}}
\makeatother
\begin{titlepage}
{~}
\vspace{-2cm}
\begin{center}
{\huge On the constants in a basic inequality for the Euler and Navier-Stokes equations}
\end{center}
\vspace{0.5truecm}
\begin{center}
{\large
Carlo Morosi$\,{}^a$, Livio Pizzocchero$\,{}^b$({\footnote{Corresponding author}})} \\
\vspace{0.5truecm} ${}^a$ Dipartimento di Matematica, Politecnico di Milano,
\\ P.za L. da Vinci 32, I-20133 Milano, Italy \\
e--mail: carlo.morosi@polimi.it \\
${}^b$ Dipartimento di Matematica, Universit\`a di Milano\\
Via C. Saldini 50, I-20133 Milano, Italy\\
and Istituto Nazionale di Fisica Nucleare, Sezione di Milano, Italy \\
e--mail: livio.pizzocchero@unimi.it
\end{center}
\begin{abstract}
We consider
the incompressible Euler or Navier-Stokes (NS) equations
on a $d$-di\-men\-sio\-nal torus $\Td$; the
quadratic term in these equations arises
from the bilinear map sending two velocity fields
$v, w : \Td \vain \reali^d$ into
$v \sc \partial w$, and also involves the
Leray projection $\LP$ onto the space of divergence free
vector fields. We derive
upper and lower bounds for the constants in some inequalities
related to the above quadratic term; these bounds hold, in particular,
for the sharp constants $K_{n d} \equiv K_n$
in the basic inequality $\| \LP(v \sc \partial w) \|_n
\leqs K_n \| v \|_n \| w \|_{n+1}$, where
$n \in (d/2, + \infty)$ and
$v, w$ are in the Sobolev spaces $\HM{n}, \HM{n+1}$
of zero mean, divergence free vector fields of orders
$n$ and $n+1$, respectively.
As examples, the numerical values of our upper
and lower bounds are reported for $d=3$ and some values of $n$.
Some practical motivations are indicated for an
accurate analysis of the constants $K_n$.
\end{abstract}
\vspace{0.2cm} \noindent
\hfill \parn
\par \vspace{0.05truecm} \noindent \textbf{AMS 2000 Subject classifications:} 76D05, 26D10, 46E35.
\end{titlepage}
\section{Introduction}
The incompressible Euler or Navier-Stokes (NS) equations in $d$ space dimensions can be written as
\beq {\partial u \over \partial t}  = - \LP(u \sc \partial u) + \nu \Delta u + f~, \label{eul} \feq
where: $u= u(x, t)$ is the divergence free velocity field; $x = (x_s)_{s=1,...,d}$ are
the space coordinates (yielding the derivatives $\partial_s := \partial/\partial x_s$);
$\Delta := \sum_{s=1}^d \partial_{s s}$ is the Laplacian;
$(u \sc \partial u)_r := \sum_{s=1}^d u_s \partial_s u_r$ ($r=1,...,d$);
$\LP$ is the Leray projection onto the space of divergence free vector fields;
$\nu = 0$ for the Euler equations; $\nu \in (0,+\infty)$ (in fact $\nu=1$,
after rescaling) for the NS equations; $f = f(x,t)$
is the Leray projected density of external forces. \parn
In this paper we stick to the case
of space periodic boundary conditions;
so, $x$ ranges in the $d$-dimensional torus
$\Td := \left({\reali /2 \pi \interi} \right)^d$.
As well known, for any solution $u$ of Eqs. \rref{eul},
the (spatial) mean $\la u \ra := (2 \pi)^{-d} \int_{\Td} u \, d x$ evolves according to
$d \la u \ra / d t = \la f \ra$,
and the zero mean vector field $u - \la u \ra$ fulfills an equation
like \rref{eul}, with $f$ replaced by a new, zero mean forcing term (see, e.g.,
\cite{due}); due to these remarks, the analysis
of Eq. \rref{eul} can be reduced to the case where
$\la u \ra = 0$, $\la f \ra = 0$.\parn
Our functional setting for the incompressible Euler/NS equations
relies on $H^n$ Sobolev spaces. More precisely we consider, for suitable (integer or noninteger)
values of $n$, the spaces
\beq \Hz{n}(\Td) \equiv \Hz{n} := \{ v : \Td \vain \reali^d~|~~
\sqrt{-\Delta}^{\,n} v \in \bb{L}^2(\Td),~ \la v \ra = 0 \}~, \feq
\beq \HM{n}(\Td) \equiv \HM{n} := \{ v \in \Hz{n}~|~\dive \, v=0 \} \feq
(the subscripts $0$, $\Sigma$ recall the
vanishing of the mean and of the divergence, respectively). For each $n$, we equip
$\Hz{n}$ with the standard inner product and the norm
\beq \la v | w \ra_n := \la \sqrt{-\Delta}^{\,n} v |  \sqrt{-\Delta}^{\,n} w \ra_{L^2}~,
\qquad \| v \|_n := \sqrt{\la v | v \ra_n}~, \feq
which can be restricted to the (closed) subspace $\HM{n}$. We can now pass to discuss
Eq. \rref{eul} with $u(\cdot, t) \in \HM{n}$ for each $t$. \parn
A fully quantitative treatment of several problems related to the above functional
setting (such as estimates on the time of existence of the solution of \rref{eul}
for a given datum, estimates on its distance from any approximate solution,
etc.) relies on the constants in some inequalities about
the bilinear map sending two vector fields
$v,w$ on $\Td$ into $v \sc \partial w$, or about the composition of
this map with $\LP$. Here, we wish to analyze the constants in some
inequality of this kind. \parn
To describe precisely the contents of this paper, let us recall that
the assumptions $n > d/2$, $v \in \HM{n}$ and $w \in \Hz{n+1}$
imply $v \sc \partial w \in \Hz{n}$, whence $\LP(v \sc \partial w) \in \HM{n}$.
In this paper we consider the basic inequality
\beq \| \LP(v \sc \partial w) \|_n \leqs K_{n} \| v \|_{n} \| w \|_{n+1}
\qquad \mbox{for $n \in (\dd{d \over 2}, + \infty)$, $v \in \HM{n}$, $w \in \HM{n+1}$}~; \label{basineq} \feq
our aim is to give quantitative
upper and lower bounds on the sharp constant
$K_{n} \equiv K_{n d}$ appearing therein. We use the fact
that $K_n \leqs K'_n$, where $K'_n$ is the sharp constant in the (auxiliary) inequality
\beq \| v \sc \partial w \|_n \leqs K'_{n} \| v \|_{n} \| w \|_{n+1}
\qquad \mbox{for $n \in (\dd{d \over 2}, + \infty)$, $v \in \HM{n}$, $w \in \Hz{n+1}$}~. \label{basinequa} \feq
Even though Eqs. \rref{basineq} \rref{basinequa} are well known,
little information can be found in the literature about the numerical values of the
constants therein. Our approach produces
fully computable upper and lower bounds $K^{\pm}_n \equiv K^{\pm}_{n d}$ such that
\beq K^{-}_n \leqs K_n \leqs K'_{n} \leqs K^{+}_n \label{uplow} \feq
for all $n > d/2$. As examples, the bounds $K^{\pm}_n$ are computed in
dimension $d=3$, for some values of $n$. In these cases the upper and lower bounds
are reasonably close, at least for the purpose to apply them to the Euler/NS equations.
\parn
In a companion paper \cite{cog}, we have proposed upper and lower bounds for the constants $G_{n d}
\equiv G_n$ in the inequality
\beq |\la v \sc \partial w | w \ra_n | \leqs G_n \| v \|_n
\| w \|^{2}_{n} \quad \mbox{for $n \in (d/2 + 1, + \infty)$, $v \in \HM{n}$, $w \in \HM{n+1}$}~,
\label{katineq} \feq
dating back to a seminal paper by Kato \cite{Kato}. \parn
Let us illustrate some applications of
the inequalities \rref{basineq} \rref{katineq}, depending on
quantitative information on the constants $K_n, G_n$.
To this purpose, following \cite{appeul} we consider the Euler/NS equations \rref{eul} with a specified
initial condition $u(x,0) = u_0(x)$; let $\ua : \Td \times [0,\Ta) \vain \reali^d$ be
an approximate solution of this Cauchy problem.
Given $n \in (d/2+1, +\infty)$ (and assuming suitable regularity
for $u_0, f, \ua$), let $\ua$ possess the differential error estimator
$\er_n : [0,\Ta) \vain [0,+\infty)$, the datum error estimator $\delta_n \in [0,+\infty)$
and the growth estimators $\Dd_n, \Dd_{n+1} : [0,\Ta) \vain [0,+\infty)$; this means that,
for $t \in [0,\Ta)$,
\beq \| \big({\partial \ua \over \partial t} +
\LP(\ua \,\sc \,\partial \ua) - \nu \Delta \ua - f \big)(t)\|_n
\leqs \er_n(t), \quad \| \ua(0) - u_0 \|_n \leqs \delta_n, \feq
$$ \| \ua(t) \|_n \leqs \Dd_n(t), \quad \| \ua(t) \|_{n+1} \leqs \Dd_{n+1}(t) $$
(with $\ua(t) := \ua(\cdot, t)$,
etc.). Furthermore, let us assume the existence of a function $\Rr_n \in C([0,\Tc),[0,+\infty))$
with $\Tc \in (0,\Ta]$, fulfilling the \textsl{control inequalities}
\beq {d^{+} \Rr_n \over d t} \geqs - \nu \Rr_n
+ (G_n \Dd_n + K_n \Dd_{n+1}) \Rr_n + G_n \Rr^2_n + \er_n
~~\mbox{on $[0,\Tc)$},~\Rr_n(0) \geqs \delta_n \label{ciq} \feq
(with $d^{+}/ d t$ the right upper Dini derivative).
Then, as shown in \cite{appeul},
the solution $u$ of Eq.\rref{eul} with initial datum $u_0$
exists (in a classical sense) on the time interval $[0,\Tc)$, and its distance from the
approximate solution admits the bound
\beq \| u(t) - \ua(t) \|_n \leqs \Rr_n(t) \qquad \mbox{for $t \in [0,\Tc)$}~. \label{buap} \feq
This somehow refines a previous result of \cite{Che}, where the time of existence
of $u$ was estimated via an integral inequality involving $\delta_n$, $\er_n$, $\Dd_n$,
$\Dd_{n+1}$ and the constants $K_n, G_n$ (but with no quantitative information
on these constants).
For a given datum $u_0$, the practical implementation of the setting
of \cite{appeul} is performed choosing a suitable $\ua$ (say, a Galerkin approximate
solution), computing the estimators $\er_n$, $\Dd_n$, $\Dd_{n+1}$ and then
using the inequalities (\ref{ciq}-\ref{buap}).
\parn
To conclude, let us mention other papers \cite{Abd,Kyr,uno,mult,due,tre,Pel,Rob,Zg1,Zg2}
where a fully quantitative approach was considered for the NS equations,
other nonlinear PDEs and/or some related inequalities, with the aim to derive
conditions of existence or error bounds on approximation methods.
In particular, in \cite{due} we derived (fairly rough) upper bounds on the constants
in a variant of the inequality \rref{basinequa} using an approach similar to the present one,
but much less refined.
\salto
\textbf{Organization of the paper.}
In Section \ref{inns} we fix our standards about Sobolev spaces on $\Td$, and introduce
in this framework the Leray projection $\LP$ and the bilinear map $v, w \mapsto v \sc \partial w$~. \parn
Section \ref{sequak} states
the main results of the paper; here we present our upper and lower bounds
$K^{\pm}_n$, fulfilling Eq.\rref{uplow} (in any space
dimension $d$); these are the subject of Propositions \ref{propup} and \ref{prolow}, respectively.
A major character of this section is a positive function $\KK_n$,
defined on the space $\Zd \setminus \{0\}$ of nonzero Fourier wave vectors,
whose sup determines our upper bound $K^{+}_n$; at each point $k \in
\Zd \setminus \{0\}$, $\KK_n(k)$ is a sum (of convolutional type)
over $\Zd \setminus \{0, k \}$. The lower bound $K^{-}_n$ given in the same section
is an elementary function of $n$ (and $d$).
As examples, the numerical values of the
bounds $K^{\pm}_{n}$ are reported for $d=3$ and $n=2,3,4,5,10$ (see Eq.
\rref{bouk}). \parn
Section \ref{provewel} contains the proofs of the previously mentioned
Propositions \ref{propup}, \ref{prolow}. \parn
Three appendices are devoted
to the practical evaluation of the function $\KK_n$ and
of the bounds $K^{+}_n$. Appendix \ref{provek}
presents some preliminary notations and results. Appendix \ref{appek} contains
the main theorem (Proposition \ref{prokknd}) about the evaluation of $\KK_n$ and of its sup.
Finally, in Appendix \ref{appe345} we give details on the computation of $\KK_n$ and $K^{+}_n$
for the previously mentioned cases $d=3$, $n=2,3,4,5,10$. \parn
For all the numerical computations required in this paper, as well as
for some lengthy symbolic manipulations, we have used systematically
the software MATHEMATICA. Throughout the paper,
an expression like $r= a. b c d e...~$ means the following: computation of
the real number $r$ via MATHEMATICA produces as an output $a.b c d e$,
followed by other digits not reported for brevity.
\section{Sobolev spaces on $\boma{\Td}$, and the Euler/NS quadratic nonlinearity}
\label{inns}
In this section we summarize our standard definitions and notations
about spaces of periodic functions and distributions, and
their applications to the incompressible Euler or NS equations;
we consider, especially, Sobolev spaces on the torus. These standards
were already described in \cite{tre}, with some more details. \parn
Throughout the paper, we consider any space dimension
\beq d \geqs 2~; \feq
$r,s$ are indices running from $1$ to $d$.
For $a = (a_r)$, $b = (b_r) \in \complessi^d$ we put
\beq a \, \sc \, b := \sum_{r=1}^d a_r \, b_r~, \qquad |a| := \sqrt{\overline{a} \, \sc \, a}~, \feq
where $\overline{a} := (\overline{a_r})$ is the complex conjugate of $a$.
Hereafter we refer to the $d$-dimensional torus
\beq \Td := \underbrace{\To \times ... \times \To}_{\tiny{\mbox{$d$ times}}}~, \qquad \To := \reali/(2 \pi \interi)~, \feq
whose elements are typically written
$x = (x_r)_{r=1,...d}$.
\salto
\textbf{Distributions on $\boma{\Td}$, Fourier series and
Sobolev spaces.}
We introduce the space of periodic distributions $\DD'(\Td, \complessi) \equiv \DD'_{\comple}$, which is
the (topological) dual of $C^{\infty}(\Td, \complessi) \equiv C^{\infty}_{\comple}$;
$\la v, f \ra \in \complessi$ denotes the action of a distribution $v \in \DD'_{\comple}$
on a test function $f \in C^{\infty}_{\comple}$. \parn
We also consider the lattice $\Zd$ of elements $k =
(k_r)_{r=1,...,d}$.
Each $v \in \DD'_{\comple}$
has a unique (weakly convergent) Fourier series expansion
\beq v = \sum_{k \in \Zd} v_k e_k~, \quad
e_k(x) := {1 \over (2 \pi)^{d/2}} \, e^{i k \sc \, x}~\mbox{for $x \in \Td$}~, \quad
v_k := \la v, e_{-k} \ra~ \in \complessi~. \label{fs} \feq
The complex conjugate of a distribution $v \in \DD'_{\comple}$ is the unique distribution $\overline{v}$ such that
$\overline{\la v, f \ra} = \la \overline{v}, \overline{f} \ra$ for each $f
\in C^{\infty}_{\comple}$; one has $\overline{v}= \sum_{k \in \Zd} \overline{v_{k}} \, e_{-k}$. \parn
The \textsl{mean} of $v \in \DD'_{\comple}$ and the space
of \textsl{zero mean} distributions are
\beq \la v \ra := {1 \over (2 \pi)^d} \la v, 1 \ra = {1 \over (2 \pi)^{d/2}} v_0~,\qquad
\DD'_{\comple \zz} := \{ v \in \DD'_{\comple}~|~\la v \ra = 0 \}
\label{mean} \feq
(of course, $\la v, 1 \ra = \int_{\Td} v \, d x$ if $v \in L^1(\Td, \complessi, d x)$).
The relevant Fourier coefficients of zero mean distributions are labeled by the set
\beq \Zd_\zz := \Zd \setminus \{0\}~. \feq
The distributional derivatives
$\partial/\partial x_s \equiv \partial_s$ and the Laplacian
$\Delta := \sum_{s=1}^d \partial_{s s}$ send $\DD'_{\comple}$ into $\DD'_{\comple \zz}$ and,
for each $v$,
$\partial_s v = i \sum_{k \in \Zd_\zz} k_s v_k e_k$,
$\Delta v = - \sum_{k \in \Zd_\zz} | k |^2 v_k e_k$.
For any $n \in \reali$, we further define
\beq \sqrt{-\Delta}^{\, n} : \DD'_{\comple} \vain \DD'_{\comple \zz}~, \qquad
v \mapsto \sqrt{- \Delta}^{\, n}  v := \sum_{k \in \Zd_\zz} | k |^{n} v_k e_k~. \label{reg} \feq
The space of
\textsl{real} distributions is
\beq \DD'(\Td, \reali) \equiv \DD' := \{ v \in \DD'_{\comple}~|~\overline{v} = v \} =
\{ v \in \DD'_{\comple}~|~\overline{v_k} = v_{-k}~\mbox{for all $k \in \Zd$} \}~. \feq
For $p \in [1,+\infty]$ we often consider the real space
\beq L^p(\Td, \reali, d x) \equiv L^p~, \qquad \feq
\parn
mainly for $p=2$. $L^2$ is a Hilbert space
with the inner product
$\la v | w \ra_{L^2} := \int_{\Td} v(x) w(x) d x = \sum_{k \in \Zd} \overline{v_k} w_k$ and
the induced norm
$\| v \|_{L^2} = \sqrt{\int_{\Td} v^2(x) d x}$ $= \sqrt{\sum_{k \in \Zd} | v_k |^2}$. \parn
The zero mean parts of $\DD'$ and $L^p$ are
\beq \DD'_\zz := \{ v \in \DD'~|~~\la v \ra = 0 \}~,
\qquad  L^p_0 := L^p \cap \Dz~; \feq
all the above mentioned differential operators send $\DD'$ into $\DD'_{\zz}$. \parn
For each $n \in \reali$, the \textsl{zero mean Sobolev space}
$H^n_\zz(\Td, \reali) \equiv H^n_\zz$ is defined by
\beq H^n_\zz := \{ v \in \DD'_\zz~|~\sqrt{-\Delta}^{\, n} v \in L^2 \}
= \{ v \in \DD'_\zz~|~\sum_{k \in \Zd_\zz} | k |^{2 n} | v_k |^2 < + \infty~\}~;
\label{hn} \feq
this is a real Hilbert space with the inner product
$\la v | w \ra_{n} := \la \sqrt{-\Delta}^{\, n} v \, | \, \sqrt{- \Delta}^{\, n} w \ra_{L^2}$
$= \sum_{k \in \Zd_\zz} | k |^{2 n} \, \overline{v_k}  w_k$ and the induced norm
$\| v \|_{n} = \| \sqrt{-\Delta}^{\, n}~ v \|_{L^2}$ $= \sqrt{\sum_{k \in \Zd_\zz} | k |^{2 n} | v_k |^2}$.
\parn
As a special case, if $n$ is a \textsl{nonnegative integer} one proves that
\beq H^n_\zz = \{ v \in \DD'_\zz~|~\partial_{s_1 ... s_n} v \in L^2~~ \forall s_1,...,s_n \in \{1,....,d\} \}
\label{hnn} \feq
and that, for each $v$ in the above space,
$\| v \|_n = \sqrt{\sum_{s_1,...,s_n=1}^d \| \partial_{s_1 ... s_n} v \|^2_{L^2}}$~.
\salto
\textbf{Spaces of vector valued functions on $\Td$.}
If $\Vi(\Td, \reali) \equiv \Vi$ is any vector space of \textsl{real}
functions or distributions on $\Td$, we write
\beq \bVi(\Td) \equiv \bVi := \{ v = (v^1,...,v^d)~|~v_r \in \Vi \quad \mbox{for all $r$}\}~. \feq
In this way we can define, e.g., the spaces $\bb{\DD}'(\Td) \equiv \bb{\DD}'$,
$\bb{L}^p(\Td) \equiv \bb{L}^p$ ($p \in [1,+\infty]$),
$\Hz{n}(\Td) \equiv \Hz{n}$.
Any $v = (v_r) \in \bb{\DD}'$ is referred to
as a (distributional) \textsl{vector field} on $\Td$. We note that $v$ has a unique Fourier
series expansion \rref{fs} with coefficients
\beq v_k := (v_{r k})_{r =1,...,d} \in \complessi^d~, \qquad v_{r k}  := \la v_r, e_{-k} \ra~; \feq
as in the scalar case, the reality of $v$ ensures $\overline{v_k} = v_{-k}$.
\parn
$\bb{L}^2$ is a real Hilbert space, with the inner product
$\la v | w \ra_{L^2} := \int_{\Td} v(x) \sc w(x) d x = \sum_{k \in \Zd} \overline{v_k} \sc w_k$
and the induced norm
\beq \| v \|_{L^2} = \sqrt{\int_{\Td} |v(x)|^2 d x} = \sqrt{\sum_{k \in \Zd} |v_k|^2}~. \label{refur} \feq
We define componentwise the mean $\la v \ra \in \reali^d$
of any $v \in \bb{\DD}'$ (see Eq. \rref{mean});
$\Dz$ is the space of zero mean vector fields, and
$\bb{L}^p_\zz = \bb{L}^p \cap \Dz$.
We similarly define componentwise the operators
$\partial_s, \Delta, \sqrt{-\Delta}^{\,n} : \bb{\DD}' \vain \Dz$. \parn
For any real $n$, the $n$-th Sobolev space of zero mean vector fields $\Hz{n}(\Td) \equiv \Hz{n}$ is made
of all $d$-uples $v$ with components
$v_r \in H^n_\zz$; an equivalent definition can be given via Eq.\rref{hn},
replacing therein $L^2$ with $\bb{L}^2$.
$\Hz{n}$ is a real Hilbert space with the inner product
$\la v | w \ra_{n} := \la \sqrt{-\Delta}^{\, n} \, v \, | \sqrt{-\Delta}^{\,n} \,w \ra_{L^2}
= \sum_{k \in \Zd_\zz} | k |^{2 n} \, \overline{v_k}  \sc \, w_k$;
the induced norm $\|~\|_{n}$ is given by
\beq \| v \|_{n} = \| \sqrt{-\Delta}^{\, n} v \|_{L^2} =
\sqrt{ \sum_{k \in \Zd_\zz} | k |^{2 n} \, |v_k|^2 }~.\feq
\salto
\textbf{Divergence free vector fields.}
Let
$\dive : \bb{\DD}' \vain \DD'_\zz$, $v \mapsto \dive \,v := \sum_{r=1}^d \partial_r v_r$
$= i \, \sum_{k \in \Zd_\zz}  (k \sc \, v_k) e_k$. Hereafter we introduce
the space $\Ds$ of \textsl{divergence free (or solenoidal) vector fields} and some
subspaces of it, putting
\beq \Ds := \{ v \in \bb{\DD}'~|~\dive \,v = 0 \}
= \{ v \in \bb{\DD}'~|~k \sc \, v_k = 0~\forall k \in \Zd~\} ~; \feq
\beq ~~\Dsz := \Ds \cap \Dz~, \quad \bb{L}^p_{\ss} := \bb{L}^p \cap \Ds~,~\bb{L}^p_{\so} := \bb{L}^p
\cap \Dsz \quad (p \in [1,+\infty])~, \feq
\beq \HM{n} := \Ds \cap \Hz{n} \quad (n \in \reali). \feq
$\HM{n}$ is a closed subspace of the Hilbert space $\Hz{n}$,
that we equip with the restrictions of  $\la~|~\ra_n$, $\|~\|_n$.
The \textsl{Leray projection} is the (surjective) map
\beq \LP : \bb{\DD}' \vain \Ds~, \qquad
v \mapsto \LP v := \sum_{k \in \Zd} (\LP_k v_k) e_k~, \label{ler1} \feq
where, for each $k$,
$\LP_k$ is the orthogonal projection of $\complessi^d$ onto the orthogonal
complement of $k$;
more explicitly, if $c \in \complessi^d$,
\beq \LP_\zz c = c~, \qquad \LP_k c = c - {k \sc \, c  \over | k |^2}\, k \quad \mbox{for $k \in \Zd_\zz$}~.
\label{ler2} \feq
From the Fourier representations of $\LP$, $\la~\ra$, etc., one easily infers
that
\beq \la \LP v \ra = \la v \ra~\mbox{for $v \in \bb{\DD}'$},
\quad \LP \Dz = \bb{D}'_{\so}, \quad \LP \bb{L}^2= \bb{L}^2_{\ss},
\quad \LP \Hz{n} = \HM{n}~\mbox{for}~ n \in \reali~. \label{tvv} \feq
Furthermore, $\LP$ is an orthogonal projection in each one of the Hilbert spaces $\bb{L}^2$, $\Hz{n}$; in
particular,
\beq \| \LP v \|_n \leqs \| v \|_n \qquad \mbox{for $v \in \Hz{n}$}~. \feq
\salto
\textbf{The quadratic Euler/NS nonlinearity.}
We are now ready to define precisely and to analyze
the bilinear map sending two (sufficiently regular) vector fields
$v, w$ on $\Td$ into $v \sc \partial w$, and the composition of this map with $\LP$.
Throughout this paragraph we assume
\beq v \in \bb{L}^2~, \quad \partial_s w \in \bb{L}^2~~ (s =1,...,d)~; \label{asin} \feq
the above condition on the derivatives of $w$ implies $w \in \bb{L}^2$. The
statements in the forthcoming two Lemmas are known, and proved only for completeness.
\begin{prop}
\label{lemma1}
\textbf{Lemma.} Consider the vector field
$v \sc \partial w$ on $\Td$, of components
\beq (v \sc \partial w)_r : = \sum_{s=1}^d v_s \partial_s w_r~; \feq
this is well defined and belongs to $\bb{L}^1$. With the additional assumption $\dive v = 0$,
one has $\la v \sc \partial w \ra = 0$ (which also implies $\la \LP(v \sc \partial w) \ra = 0$, see
\rref{tvv}).
\end{prop}
\textbf{Proof.} Each component $(v \sc \partial w)_r$, being a sum of
products of $L^2$ functions, is evidently in $L^1$. \parn
Now, assume provisionally that $v, w$ are $C^1$; then
$\int_{\Td} (v \sc \partial w)_r \, d x$ $= \sum_{s=1}^d \int_{\Td}  v_s \partial_s w_r \, d x $ $=
- \sum_{s=1}^d \int_{\Td} (\partial_s v_s)  w_r \, d x $  (integrating by parts), i.e.,
\beq \int_{\Td} (v \sc \partial w)_r \, d x = - \int_{\Td}  (\dive v) \, w_r \, d x~. \label{densit} \feq
By simple density arguments, \rref{densit} holds whenever $v, \partial_s w \in \bb{L}^2$
and $\dive v \in L^2$. In particular, $\int_{\Td} v \sc \partial w \, d x = 0$ if
$v, \partial_s w \in \bb{L}^2$ and $\dive v=0$. \fine
\begin{prop}
\label{lemma3}
\textbf{Lemma.}
$v \sc \partial w$ has Fourier coefficients
\beq (v \sc \, \partial w)_k =  {i \over (2 \pi)^{d/2}}
\sum_{h \in \Zd} [v_{h} \sc \, (k - h)] w_{k - h}
\qquad \mbox{for all $k \in \Zd$}~. \label{dainf} \feq
\end{prop}
\textbf{Proof} (Sketch).
Consider the Fourier coefficients $v_{r k}, w_{s k}$; then $(\partial_r w_s)_k =
i k_r \, w_{s k}$. The pointwise product corresponds to $(2 \pi)^{-d/2}$ times the convolution
of the Fourier coefficients; thus
$(v \sc \, \partial w)_{s k}$ $=  i (2 \pi)^{-d/2} \sum_{r=1}^d$
$\sum_{h \in \Zd} v_{r h} (k - h)_r w_{s, k - h}$; the
vector form of this statement is Eq. \rref{dainf}. \fine
To conclude, we note that Eqs. \rref{ler1} \rref{dainf} imply
\beq [\LP(v \sc \partial w)]_k = {i \over (2 \pi)^{d/2}} \sum_{h \in \Zd} [v_{h} \sc \, (k - h)] \LP_k w_{k - h}
\qquad \mbox{for all $k \in \Zd$}~, \label{infert} \feq
with $\LP_k$ as in \rref{ler2}.
\section{The basic inequality for the Euler/NS quadratic nonlinearity}
\label{sequak}
Throughout this section we assume ($d \in \{2,3,...\}$ and)
\beq n \in ({d \over 2}, + \infty)~. \label{assd} \feq
Given two vector fields $v, w$ on $\Td$, we have already discussed
$v \sc \partial w$
under the conditions $v, \partial_s w \in \bb{L}^2$; here we consider the much stronger
assumptions $v$ in $\HM{n}$, $w$ in $\Hz{n+1}$ or $\HM{n+1}$. \parn
The forthcoming Proposition \ref{weldef} is well known, and
presented here only for completeness; as a matter of fact,
the quantitative analysis performed later
will also give, as a byproduct, an alternative proof of this Proposition.
\begin{prop}
\label{weldef}
\textbf{Proposition.}
Let $v \in \HM{n}$, $w \in \Hz{n+1}$; then,
\beq  v \sc \partial w \in \Hz{n}~. \feq
The map
\beq \HM{n} \times \Hz{n+1} \vain \Hz{n}, \qquad (v, w) \mapsto
v \sc  \partial w  \feq
is bilinear and continuous.
\end{prop}
Of course, continuity of the above map
is equivalent to the existence of a
nonnegative constant $K'$, such that $\| v \sc  \partial w \|_n \leqs K' \| v \|_n
\| w \|_{n+1}$ for $v, w$ as in the previous Proposition. A similar inequality
holds as well for $\LP(v \sc \partial w) \in \HM{n}$, since $\| \LP(v \sc  \partial w) \|_n
\leqs \| v \sc  \partial w \|_n$. \parn
So, we have the ''auxiliary inequality'' \rref{basinequa}
and the ''basic inequality'' \rref{basineq} of the Introduction; the sharp constants
appearing therein can be defined as follows.
\begin{prop}
\textbf{Definition.} We put
\beq K'_{n d} \equiv K'_n  \label{kpnd} \feq
$$ := \min \{ K' \in [0,+\infty)~|~\| v \sc \partial w \|_n \leqs K' \| v \|_n
\| w \|_{n+1}~\mbox{for all $v \in \HM{n}$, $w \in \Hz{n+1}$}\}~; $$
\beq K_{n d} \equiv K_n  \label{knd} \feq
$$ := \min \{ K \in [0,+\infty)~|~\| \LP(v \sc \partial w) \|_n \leqs K \| v \|_n
\| w \|_{n+1}~\mbox{for all $v \in \HM{n}$, $w \in \HM{n+1}$}\}~. $$
(Note that all $w$'s in \rref{knd} are divergence free, a property not required in
\rref{kpnd}.)
\end{prop}
The considerations after Proposition \ref{weldef} ensure that
\beq K_n \leqs K'_n \,; \label{ensure} \feq
in the rest of the section (which is its original part) we present
computable upper and lower bounds on $K'_n$ and $K_{n}$, respectively. \parn
The upper bound requires a more lengthy analysis; the final result relies on a
function $\KK_{n d} \equiv \KK_{n}$, appearing in the forthcoming Definition \ref{dekknd}.
Hereafter we introduce some auxiliary notations, required to build $\KK_n$.
\begin{prop}
\label{area}
\textbf{Definition.}
(i) Here and in the sequel, the exterior power $\We^2 \reali^d$
is identified with the space of real, skew-symmetric $d \times d$ matrices
$A = (A_{r s})_{r,s,=1,...,d}$; this is equipped with the operation of exterior product
\beq \we : \reali^d \times \reali^d \vain {\We}^2 \reali^d,
\quad (p,  q) \mapsto p \we q~~\mbox{s.t.}~~
(p \we q)_{r s} := p_r q_s - q_r p_s~.
\label{ester} \feq
(ii) We equip the above space with the norm
\beq |~| : {\We}^2 \reali^d \vain [0,+\infty),
\qquad A = (A_{r s}) \mapsto |A| :=  \sqrt{{1 \over 2}\, \sum_{r,s=1}^d  |A_{r s}|^2}~.
\label{normest}
\feq
\end{prop}
The operation \rref{ester} is bilinear and skew-symmetric; when composed with the norm
\rref{normest}, it gives a mapping
\beq \reali^d \times \reali^d \vain [0,+\infty)~,
\qquad (p,q) \mapsto |p \we q|~, \label{sqp} \feq
which has the following, well known properties.
\begin{prop}
\textbf{Proposition.}
Let $p, q$ in $\reali^d$. Then
\beq |p \we q| = \sqrt{|p|^2 |q|^2 - (p \sc q)^2}~=
|p| |q| \sin \tet~,
\label{norqp} \feq
where $\tet \equiv \tet(p, q) \in [0; \pi]$ is the convex angle between $p$ and $q$ (defined arbitrarily, if
$p = 0$ or $q = 0$). So,
$|p \we q|$ is the area of the parallelogram of sides $p$,$q$ and
\beq |p \we q|  \leqs |p||q|~. \feq
\end{prop}
Now we are ready to construct the function $\KK_{n}$, a major character of
the section.
\begin{prop}
\label{dekknd}
\textbf{Definition.} We put
\beq \Zd_{0 k} := \Zd \setminus \{0, k \} \qquad \mbox{for each $k \in \Zd_0$}~;  \label{zdok} \feq
\beq \KK_{n d} \equiv \KK_n : \Zd_0 \vain (0,+\infty),\quad
k \mapsto  \KK_{n}(k) := |k|^{2 n} \sum_{h \in \Zd_{0 k}} {|h \we k|^2 \over |h|^{2 n + 2} |k-h|^{2 n + 2}}.
\label{kknd} \feq
\end{prop}
\begin{rema}
\textbf{Remarks.}
(i) The sum in \rref{kknd} has nonnegative terms, so it exists
in principle as an element of $[0,+\infty]$. However,
\beq \sum_{h \in \Zd_{0 k}} {|h \we k|^2 \over |h|^{2 n + 2} |k-h|^{2 n + 2}} < + \infty~
\mbox{for all $k \in \Zd_0$}~; \label{conve} \feq
in fact, for fixed $k$ and $h \vain \infty$,
\beq {|h \we k|^2 \over |h|^{2 n + 2} |k-h|^{2 n + 2}} = O\Big({1 \over |h|^{4 n + 2}}\Big)
\feq
and, for each family $(s_h)_{h \in \Zd}$ with elements
in $[0,+\infty)$, the relation $s_h = O(1/|h|^\nu)$ with $\nu > d$ implies,
as well known, $\sum_{h \in \Zd} s_h < + \infty$. These considerations justify the claim
$\KK_{n}(k) \in (0,+\infty)$ in \rref{kknd}. \parn
(ii) In Eq. \rref{kknd}, one can insert at will the identity
$h \we k = h \we (k -h)$
(following from the bilinearity and skew-symmetry of $\we$), and
the inequality $| h \we (k - h) | \leqs |h | |k-h|$;
this will be occasionally done in the sequel. \parn
(iii) Let $r \in \{1,...,d\}$, and let $\sigma$ be any permutation of $\{1,...,d\}$;
introduce the reflection operator $R_r$ and the permutation operator $P_\si$
defined by
\beq R_r, P_\si : \reali^d \vain \reali^d ~, \label{refper} \feq
$$ R_r(k_1,.., k_r,...,k_d) := (k_1,...,-k_r,...,k_d)~, \qquad P_{\si}(k_1,...,k_d) :=
(k_{\si(1)},...,k_{\si(d)})~; $$
these are orthogonal operators (with respect to the inner product $\sc$ of $\reali^d$),
and send $\Zd_0$ into itself. We note that
\beq \KK_{n}(R_r k) = \KK_{n}(k)~, \qquad \KK_{n}(P_\sigma k) = \KK_{n}(k)~
\qquad \mbox{for each $k \in \Zd_0$}~; \label{claim} \feq
for example, the first equality is checked expressing
$\KK_{n}(R_r k)$ via the definition \rref{kknd}, making a change of
variable $h \mapsto R_r h$ in the sum therein and noting that $|(R_r h) \we (R_r k)| =
|h \we k|$, $|R_r h| = |h|$, $|R_r k|= |k|$, $|R_r k - R_r h| = |k-h|$.
The verification of the second inequality \rref{claim} proceeds similarly. \parn
Due to the symmetry properties \rref{claim}, the computation of $\KK_{n}(k)$
can always be reduced to the case $k_1 \geqs k_2 \geqs ... \geqs k_d \geqs 0$.
\parn
(iv) In Appendix \ref{appek} we prove that
\beq \sup_{k \in \Zd_0} \KK_{n}(k) < + \infty~. \label{finito} \feq
This appendix also gives tools for the practical evaluation of $\KK_{n}$ and of its sup.
\fine
\end{rema}
Let us pass to the desired upper bound, which is the following.
\begin{prop}
\label{propup}
\textbf{Proposition.}
The constant $K'_{n}$ defined by \rref{kpnd} has the upper bound
\beq K'_{n} \leqs K^{+}_n~, \label{hasup} \feq
\beq K^{+}_{n} := {1 \over (2 \pi)^{d/2}}
\sqrt{\sup_{k \in \Zd_0} \KK_{n}(k)}~~(\mbox{or any upper approximant for this}). \label{deup} \feq
\end{prop}
\textbf{Proof.} See Section \ref{provewel}. \fine
Let us pass to the problem of finding a
lower bound for $K_{n}$; this can be obtained
directly from the tautological inequality $K_n \geqs \| v \sc \partial w \|_n /\| v \|_n \| w \|_{n+1}$,
choosing for $v, w$ some suitable trial functions. A very simple choice of $v, w$ yields
the following.
\begin{prop}
\label{prolow}
\textbf{Proposition.} The constant $K_{n}$ defined by \rref{knd} has the lower bound
\beq K_{n} \geqs K^{-}_n~, \label{resulow} \feq
\beq K^{-}_n := {2^{n/2} \over (2 \pi)^{d/2}} U_d ~~(\mbox{or any round down for this number})~, \label{delow}
\feq
$$ U_d := \left\{\barray{cc} (2 - \sqrt{2})^{1/2} = 0.76536...
& \mbox{if $d=2$,} \\ 1 & \mbox{if $d \geqs 3$.} \farray \right. $$
\end{prop}
\textbf{Proof.} See Section \ref{provewel}. \fine
Putting together Eqs. \rref{ensure} \rref{hasup} \rref{resulow} we obtain
a chain of inequalities, anticipated in the Introduction,
$$ K^{-}_n \leqs K_n \leqs K'_{n} \leqs K^{+}_n~; $$
here, the bounds $K^{\pm}_n$ can be computed explicitly from their definitions
\rref{deup} \rref{delow}.
\begin{rema}
\textbf{Examples.}
For $d=3$ and $n=2,3,4,5,10$, we can take
\beq K^{-}_{2} = 0.126~,~~K^{+}_2 = 0.335~; \qquad
K^{-}_{3} = 0.179~,~~ K^{+}_{3} = 0.323~, \label{bouk} \feq
$$ K^{-}_{4} =0.253~,~~K^{+}_{4} =  0.441~; \quad
K^{-}_5 = 0.359~,~~K^{+}_{5} = 0.657~; \quad K^{-}_{10} = 2.03~,~~ K^{+}_{10} = 6.21~. $$
In the above, the $K^{-}_n$ are obtained rounding down to three digits the
number ${2^{n/2} (2 \pi)^{-3/2}}$; the $K^{+}_n$ are obtained from upper approximation
of the sup in \rref{deup}, as illustrated in Appendix \ref{appe345}.
The ratios $K^{-}_{n}/K^{+}_n$
are  $0.376...$, $0.554...$, $0.573...$, $0.546...$,
$0.326...$ for $n=2,3,4,5,10$, respectively. One can see that
$K^{-}_{n}/K^{+}_n$ is smaller (i.e., that we have a larger
uncertainty on the sharp constant $K_n$) in the extreme cases $n=2$, $n=10$; we presume
this to happen, in any space dimension $d$, when $n$ approaches
the limit values $d/2$ and $+\infty$.
\parn
To avoid misunderstandings related to the above examples,
we repeat that the approach of this paper
applies as well to noninteger values of $n$.
\end{rema}
\section{Proof of Propositions (\ref{weldef} and) \ref{propup}, \ref{prolow}}
\label{provewel}
Throughout the section $n \in (\dd{d/2},+\infty)$.
\begin{prop}
\label{lemabz}
\textbf{Lemma.} Let
\beq p, q \in \reali^d \setminus \{0 \}~, \quad z \in \complessi^d,~p \sc z = 0~, \feq
and $\tet(p,q) \equiv \tet \in [0,\pi] $ be the convex angle between $q$ and $p$. Then
\beq | q \sc z | \leqs \sin \tet \,  |q | |z | =
{|p \we q| \over |p|} \, |z |~.  \label{qscz} \feq
\end{prop}
\textbf{Proof.}
We choose an orthonormal basis $(\eta_r)_{r=1,...,d}$ of $\reali^d$
so that $q$ be a positive multiple of $\eta_1$, $p$ be in the span of $\eta_1, \eta_2$ and
$p \sc \eta_2 \geqs 0$; then
\beq q = |q| \eta_1~,
\qquad p = |p| (\cos \tet  \, \eta_1 + \sin \tet \, \eta_2)~. \label{espb} \feq
The $(d-1)$ vectors
\beq - \sin \tet \, \eta_1 + \cos \tet \, \eta_2, \eta_3, ... , \eta_d \feq
clearly form an orthonormal  basis for
\beq \{ p \}^{\perp} := \{ z \in \complessi^d~|~p \sc z = 0 \}~; \feq
so, any $z \in \{ p \}^{\perp}$ has a unique expansion
\beq z = z^{(2)} (- \sin \tet  \, \eta_1 + \cos \tet \, \eta_2) + z^{(3)} \eta_3 + ...
+ z^{(d)} \eta_{d},~z^{(t)} \in \complessi~\mbox{for $t=2,...,d$}~.  \label{espz} \feq
From Eqs. \rref{espb} \rref{espz} we get
\beq q \sc z = - \sin \tet \,|q| z^{(2)}~, \feq
which implies
\beq | q \sc z | = \sin \tet \,|q| |z^{(2)}| \leqs \sin \tet |q| |z|~. \feq
So, the inequality in \rref{qscz} is proved;
the subsequent equality in \rref{qscz} follows from \rref{norqp}.
\fine
\noindent
\textbf{Proofs of Propositions \ref{weldef} and \ref{propup}.} We choose
$v \in \HM{n}$, $w \in \Hz{n+1}$ and proceed in two steps; let us
recall that $v \sc \partial w$ has zero mean, see Lemma \ref{lemma1}. \parn
\textsl{Step 1. The Fourier coefficients of $v \sc \partial w$, and some
estimates for them.} First of all $(v \sc \partial w)_0=0$. The other Fourier coefficients are
\beq (v \sc \partial w)_k = {i \over (2 \pi)^{d/2}} \sum_{h \in \Zd_{0 k}} [v_{h} \sc \, (k - h)] w_{k - h}
\qquad \mbox{for $k \in \Zd_0$}~; \label{infertt} \feq
this follows from \rref{infert} taking into account that, in the sum therein,
the term with $h=0$ vanishes due to $v_0=0$, and the term with $h = k$ is zero for evident reasons.
\parn
Taking \rref{infertt} as a starting point, let us
make some remarks on the term $v_{h} \sc \, (k - h)$ appearing
therein. We have $h \sc v_h = 0$ due to the assumption $\dive v=0$;
so, we can apply Eq. \rref{qscz} with $p=h$, $q=k-h$ and $z = v_h$, which gives
\beq |v_{h} \sc \, (k - h)| \leqs {|h \we (k-h)| \over |h|} |v_h| =
{|h \we k| \over |h|} |v_h| \label{ins1} \feq
(recall that $h \we (k-h) = h \we k$). \parn
Eqs. \rref{infertt} and \rref{ins1} imply the following, for each $k \in \Zd_0$:
\parn
\vbox{
$$ | (v \sc \partial w)_k | \leqs {1 \over (2 \pi)^{d/2}} \sum_{h \in \Zd_{0 k}}
{|h \we k| \over |h|} | v_h | | w_{k - h} | $$
\beq = {1 \over (2 \pi)^{d/2}} \sum_{h \in \Zd_{0 k}}
{|h \we k| \over |h|^{n + 1} |k - h |^{n + 1}} \Big( |h |^n | v_h | |k - h|^{n+1} | w_{k - h} | \Big)~;
\label{ins4} \feq
}
\noindent
now, H\"older's inequality $| \sum_h~ a_h b_h |^2 \leqs \Big(\sum_h | a_h |^2\Big)
\Big(\sum_h~| b_h |^2 \Big)$ gives
\parn
\vbox{
\beq | (v \sc \partial w)_k |^2 \leqs {1 \over (2 \pi)^d} \, \CA_{n}(k) \CB_{n}(k)~\mbox{for
all $k \in \Zd_0$}, \label{dains} \feq
$$ \CA_{n}(k) :=  \sum_{h \in \Zd_{0 k}} {|h \we k|^2 \over | h |^{2 n + 2} | k - h |^{2 n + 2} }~, $$
$$ \CB_{n}(k) \equiv \CB_{n}(v,w)(k) := \sum_{h \in \Zd_{0 k}} |h|^{2 n} |  v_{h} |^2
| k - h |^{2 n+2} | w_{k - h} |^2 $$
}
\noindent
(in the definition of $\CB_{n}(k)$ one can write as well $\sum_{h \in \Zd_{0}}$, since
the general term of the sum vanishes for $h=k$).
We now multiply both sides of \rref{dains} by $|k|^{2 n}$; it appears that $|k|^{2 n} \CA_{n}(k) =
\KK_{n}(k)$ with $\KK_{n}(k)$ as in \rref{kknd}, so
\beq | k |^{2 n} |(v \sc \partial w)_k |^2 \leqs {1 \over (2 \pi)^d} \KK_{n}(k) \CB_{n}(k)~. \label{dainsk} \feq
\textsl{Step 2. Completing the proofs of Propositions \ref{weldef}, \ref{propup}}.
Due to \rref{dainsk},
$$ \sum_{k \in \Zd_0} | k |^{2 n} |(v \sc \partial w)_k |^2 \leqs
{1 \over (2 \pi)^d} \sum_{k \in \Zd_0} \KK_{n}(k) \CB_{n}(k) \leqs
{1 \over (2 \pi)^d} \Big(\sup_{k \in \Zd_0} \KK_{n}(k)\Big) \Big(\sum_{k \in \Zd_0} \CB_{n}(k) \Big)~. $$
The sup of $\KK_n$ is finite, as we will show (by an independent
argument) in Proposition \ref{prokknd}.
Making reference to the definition of $K^{+}_n$ in terms of this sup (see Eq. \rref{deup}),
we can write the last result as
\beq \sum_{k \in \Zd_0} | k |^{2 n} |(v \sc \partial w)_k |^2
\leqs (K^{+}_n)^2 \sum_{k \in \Zd_0} \CB_{n}(k)~. \label{daret} \feq
On the other hand, making explicit the definition of
$\CB_{n}$ we see that
\beq \sum_{k \in \Zd_0} \CB_{n}(k) = \sum_{k \in \Zd_0} \sum_{h \in \Zd_0} |h|^{2 n} |  v_{h} |^2
| k - h |^{2 n+2} | w_{k - h} |^2 \feq
$$ = \sum_{h \in \Zd_0} |h|^{2 n} |  v_{h} |^2 \sum_{k \in \Zd_{0}} | k - h |^{2 (n+1)} | w_{k - h} |^2
= \Big(\sum_{h \in \Zd_0} |h|^{2 n} |  v_{h} |^2 \Big)
\Big(\sum_{\ell \in \Zd_{0 h}} | \ell |^{2 (n+1)} | w_{\ell} |^2 \Big) $$
$$  \leqs \Big(\sum_{h \in \Zd_0} |h|^{2 n} |  v_{h} |^2 \Big)~
\Big(\sum_{\ell \in \Zd_{0}} | \ell |^{2 (n+1)} | w_{\ell} |^2 \Big) = \| v \|^2_{n} \, \| w \|^2_{n+1}~. $$
Returning to \rref{daret}, we obtain
\beq \sum_{k \in \Zd_0} | k |^{2 n} |(v \sc \partial w)_k |^2 \leqs (K^{+}_n)^2
\| v \|^2_{n} \, \| w \|^2_{n+1}~. \label{tells} \feq
We already know that $v \sc \partial w$ has zero mean. Eq.
\rref{tells} indicates the finiteness of $\sum_{k \in \Zd_0} | k |^{2 n} |(v \sc \partial w)_k |^2$, so
\beq v \sc \partial w \in \Hz{n}~; \feq
Eq. \rref{tells} also gives
\beq \| v \sc \partial w \|_n \leqs K^{+}_n \| v \|_n \| w \|_{n+1}~. \label{algives} \feq
Now, we let $(v,w)$ vary. The map
$\HM{n} \times \Hz{n+1} \vain \Hz{n}$, $(v, w) \mapsto v \sc \partial w$
is clearly bilinear, and \rref{algives} indicates its continuity; so, Proposition \ref{weldef} is proved. \parn
Eq. \rref{algives} indicates as well that the sharp constant $K'_n$ in the inequality
$\| v \sc \partial w \|_{n} \leqs K'_{n} \| v \|_n \| w \|_{n+1}$ fulfills
$K'_n \leqs K^{+}_n$, thus proving Eq. \rref{hasup} and Proposition \ref{propup}. \fine
\textbf{Proof of Proposition \ref{prolow}.}
Consider any $v \in \HM{n} \setminus \{0 \}, w \in \HM{n+1} \setminus \{0 \}$; then
\beq K_{n} \geqs {\| \LP(v \sc \partial w) \|_n \over \| v \|_n \| w \|_{n+1}}~. \label{eqthen} \feq
Hereafter we choose $v, w$ with Fourier coefficients
\beq v_k = A \delta_{k, a} + \Ac \delta_{k, -a}~, \quad w_k = B \delta_{k, b} + \Bc \delta_{k, -b} \label{defuv} \feq
($\delta$ the usual Kronecker symbol), where
\beq a := (1,0,...,0), \quad b := (0,1,0,...,0)  \label{defuv1} \feq
\beq A, B \in \complessi^d \setminus \{0\}~,~~A \sc a = 0~, \quad B \sc b = 0~; \label{defuv2} \feq
the above conditions on $A, B$ are fulfilled if and only if
\beq A = (0,\al,\dal),~B = (\be, 0,\dbe), \quad
\al, \be \in \complessi,~\dal, \dbe \in \complessi^{d-2},~ (\al,\dal), (\be,\dbe)
\neq (0,0)~. \label{defuv3} \feq
(In the case $d=2$, one understands $\dal, \dbe$ to be missing from the
above formulas: $A= (0,\al)$, $B=(\be,0$).) \parn
Of course, Eqs. (\ref{defuv}-\ref{defuv2}) ensure $v,w \in \HM{m}$ for all $m \in \naturali$
(in particular, $v,w$ are divergence free due to $a\sc A=0$, $b\sc B=0$). We now compute
the right hand side of Eq. \rref{eqthen}, in several steps. \parn
\textsl{Step 1. The norms $\| v \|_n$, $\| w \|_{n+1}$}. From the Fourier representation of
$\|~\|_n$ and from \rref{defuv}, one gets
$ \| v \|^2_n = | a |^{2 n} |A |^2 + |a |^{2 n} |\Ac|^2$, whence (noting that $|a|=1$)
\beq  \| v \|^2_n = 2 |A|^2~= 2 (|\al|^2 + |\dal|^2)~; \label{norv} \feq
similarly,
\beq  \| w \|^2_{n+1} = 2 |B|^2 = 2(|\be|^2 + |\dbe|^2)~. \label{norw} \feq
\textsl{Step 2. The Fourier coefficients of $\LP(v \sc \partial w)$}.
Let $k \in \Zd_0$; from Eqs. \rref{infert} and
(\ref{defuv}-\ref{defuv2}) we get
\beq (2 \pi)^{d/2} [\LP(v \sc \partial w)]_k = i \sum_{h = \pm a}  [v_{h} \sc \, (k - h)]  \LP_k w_{k - h}   \label{weget} \feq
$$ =  i [A \sc \, (k - a)] \LP_k w_{k - a} + i [\Ac \sc \, (k + a)] \LP_k w_{k + a} =
 i (A \sc \, k) \LP_k w_{k - a} + i (\Ac \sc \, k) \LP_k w_{k + a}  $$
$$ =  i (A \sc \, k) (\delta_{k-a, b} \LP_k B + \delta_{k-a, -b} \LP_k \Bc)
 + i (\Ac \sc \, k) (\delta_{k+a, b} \LP_k B + \delta_{k+a, -b} \LP_k \Bc)~.$$
On the other hand,
\parn
\vbox{
$$ (A \sc \, k) \delta_{k-a, b} \LP_k B =  (A \sc b) \delta_{k, a + b} \LP_{a + b} B~,$$
$$ (A \sc \, k) \delta_{k-a, -b} \LP_k \Bc = -(A \sc \, b) \delta_{k,a -b} \LP_{a-b} \Bc~, $$
$$ (\Ac \sc \, k) \delta_{k+a, b} \LP_k B = (\Ac \sc b) \delta_{k,-a + b} \LP_{-a +b} B =
(\Ac \sc b) \delta_{k,-a + b} \LP_{a -b} B~, $$
$$  (\Ac \sc \, k) \delta_{k+a, -b} \LP_k \Bc
= -(\Ac \sc \, b) \delta_{k,-a -b} \LP_{-a-b} \Bc
= -(\Ac \sc \, b) \delta_{k,-a - b} \LP_{a+b} \Bc~;
$$}
\noindent
so, returning to \rref{weget} we get
\beq (2 \pi)^{d/2} [\LP(v \sc \partial w)]_k    \label{gett} \feq
$$ = i (A \sc \, b) (\delta_{k, a + b} \LP_{a + b} B - \delta_{k,a -b} \LP_{a-b} \Bc)
 - i (\Ac \sc \, b) ( \delta_{k,-a -b} \LP_{a+b} \Bc -\delta_{k,-a + b} \LP_{a - b} B)~. $$
The explicit expressions \rref{defuv1} \rref{defuv2} for $a, b, A, B$ give
\beq A \sc b = \al,~~\Ac \sc b = \alc~; \label{ascb} \feq
\beq \LP_{a \pm b} B = B - {(a \pm b) \sc B \over |a  \pm b|}(a \pm b)
= (\be,0,\dbe) - {\be \over \sqrt{2}} (1,\pm 1,0,...,0) \label{equaa} \feq
$$= \left(\left(1 - {1 \over \sqrt{2}}\right) \be, \mp {\be \over \sqrt{2}}, \dbe\right)~; $$
\beq \LP_{a \pm b} \Bc = \overline{\LP_{a \pm b} B} =
\left(\left(1 - {1 \over \sqrt{2}}\right) \bec, \mp {\bec \over \sqrt{2}}, \dbec\right)~. \label{equab} \feq
Inserting Eqs. (\ref{ascb}-\ref{equab}) into \rref{gett}, we obtain the final result
$$ [\LP(v \sc \partial w)]_k  = {i \al \over (2 \pi)^{d/2}} \left[
\Big(\Big(1 - {1 \over \sqrt{2}}\Big) \be, - {\be \over \sqrt{2}}, \dbe\Big)
\delta_{k, a + b}
- \Big(\Big(1 - {1 \over \sqrt{2}}\Big) \bec, {\bec \over \sqrt{2}}, \dbec\Big)
\delta_{k, a - b} \right] $$
\beq - {i \alc \over (2 \pi)^{d/2}} \left[
\Big(\Big(1 - {1 \over \sqrt{2}}\Big) \bec, - {\bec \over \sqrt{2}}, \dbec\Big)
\delta_{k, - a - b}
- \Big(\Big(1 - {1 \over \sqrt{2}}\Big) \be, {\be \over \sqrt{2}}, \dbe\Big)
\delta_{k,-a+b} \right]. \label{finares} \feq
\textsl{Step 3. The norm of $\LP(v \sc \partial w)$.} From \rref{finares} and the Fourier representation
of $\|~\|_n$ we get
\parn
\vbox{
\beq \| \LP(v \sc \partial w) \|^{2}_n = \sum_{k=a \pm b,-a \mp b} |k|^{2 n}
|[\LP(v \sc \partial w)]_k|^2 \label{idest} \feq
$$ = {2^n \over (2 \pi)^d} ~ 4 |\al|^2
\Big( \Big(1 - {1 \over \sqrt{2}}\Big)^2 |\be|^2 + {1 \over 2}|\be|^2 + |\dbe|^2 \Big)
= {2^{n+2} \over (2 \pi)^d} |\al|^2 \left( (2 - \sqrt{2}) |\be|^2 + |\dbe|^2 \right)~.$$
}
\noindent
\textsl{Step 4. The lower bound on $K_{n}$.} We return to the inequality \rref{eqthen}, using
the expressions \rref{norv}, \rref{norw}, \rref{idest} for the norms of $v,w, \LP(v \sc \partial w)$; this gives
\beq K^{2}_{n} \geqs {2^{n} \over (2 \pi)^d}
{|\al|^2 \left( (2 - \sqrt{2}) |\be|^2 + |\dbe|^2 \right) \over (|\al|^2 + |\dal|^2) (|\be|^2 + |\dbe|^2)}
\label{lowb} \feq
for all $(\al,\dal), (\be, \dbe) \in \complessi \times \complessi^{d-2} \setminus \{(0,0)\}$. \parn
In the case $d=2$, one
understands $\dal,\dbe$ to be missing from the above formula; so, \rref{lowb} gives
\beq  K^{2}_{n} \geqs {2^{n} \over (2 \pi)^2} (2 - \sqrt{2})~. \label{resu2} \feq
In the case $d \geqs 3$, we choose $(\al,\dal), (\be, \dbe) \neq (0,0)$ so as to maximize the
right hand side of Eq. \rref{lowb}. The maximum is attained with $\dal = 0$, $\be=0$ and arbitrary
$\al \in \complessi \setminus \{0\}$, $\dbe \in \complessi^{d-2} \setminus \{0\}$; this choice gives
\beq K^{2}_{n} \geqs {2^{n} \over (2 \pi)^d} \qquad (d \geqs 3). \label{resu3} \feq
The results \rref{resu2} \rref{resu3} are summarized by Eqs. (\ref{resulow}-\ref{delow})
in the statement of the Proposition, which is now proved. \fine
\noindent
\appendix
\section{Some tools preparing the analysis of the function $\boma{\KK_{n}}$}
\label{provek}
Let us fix some notations, to be used throughout the Appendices.
\begin{prop}
\textbf{Definition.}
\label{versk}
(i) $\teta : \reali \vain \{0,1\}$
is the Heaviside function such that $\teta(z) := 1$ if $z \in [0,+\infty)$ and $\teta(z) := 0$
if $z \in (-\infty,0)$. \parn
(ii) $\Gamma$ is the Euler Gamma function,
$\binom{\cdot}{\cdot}$ are the binomial coefficients. \parn
(iii) We put $\Sd := \{ u \in \reali^d~|~|u| = 1 \}$.
For each $p \in \reali^d \setminus \{0\}$, the versor of $p$ is
$\vers{p}  := \dd{p \over |p|} \in \Sd$.
\end{prop}
\begin{prop}
\label{lemmaf}
\textbf{Lemma.} For any function $f : \Zd_0 \vain \reali$ and $k \in \Zd_0$,
$\ro \in (1,+\infty)$, one has
\beq \sum_{h \in \Zd_{0 k}, |h| < \ro \op |k-h| < \ro} f(h) =
\sum_{h \in \Zd_{0 k},|h| < \ro} f(h) + \teta(|k-h| - \ro) f(k-h)~.
\label{tesif} \feq
\end{prop}
\textbf{Proof.} The domain of the sum is the disjoint union of the sets
$\{ h \in \Zd_{0 k}~|~|h| < \ro \}$ and
$\{ h \in \Zd_{0 k}~|~|k-h| < \ro, |h| \geqs \ro \}$;
so,
$$ \sum_{h \in \Zd_0, |h| < \ro \op |k-h| < \ro} f(h) = \sum_{h \in \Zd_{0 k}, |h| < \ro} \!\!\!\!\!\!f(h) +
\sum_{h \in \Zd_{0 k}, |k-h| < \ro, |h| \geqs \ro}  \!\!\!f(h)  $$
$$ = \sum_{h \in \Zd_{0 k}, |h| < \ro} \!\!\!f(h) +
\sum_{h \in \Zd_{0 k}, |k-h| < \ro} \!\!\!\teta(|h| - \ro)\, f(h) ~.
$$
Now, a change of variable $h \vain k - h$ in the last sum gives
the thesis \rref{tesif}. \fine
\begin{prop}
\label{lembn}
\textbf{Lemma.}
For any $n \in [1,+\infty)$ and $p, q \in \reali^d$, one has
\beq |p \we q|^2 | p + q |^{2 n} \leqs {2^{2 n+1} (n+1)^{n+1} \over (n+2)^{n+2}}~
|p|^2 |q|^2 ( |p|^{2 n} + |q|^{2 n} )~.
\label{onthe} \feq
\end{prop}
\textbf{Proof.}
Eq. \rref{onthe} is obvious if $p=0$ or $q=0$, due to the vanishing of both sides;
hereafter we prove \rref{onthe} for $p, q \in \reali^d \setminus \{0\}$.
Let $\tet \in [0,\pi]$ denote the convex angle between
$p$ and $q$; then
$$ |p \we q|^2 = |p|^2 |q|^2 (1-\cos^2 \tet)~, \quad
|p + q|^2 = |p|^2 + |q|^2 + 2 |p| |q| \cos \tet~, $$
so
\parn
\vbox{
\beq {|p \we q|^2 | p + q |^{2 n} \over |p|^2 |q|^2 ( |p|^{2 n} + |q|^{2 n} )}
\label{eqgivea} \feq
$$ = {(1 - \cos^2 \tet) (|p|^2 + |q|^2 + 2 |p| |q| \cos \tet)^n
\over  |p|^{2 n} + |q|^{2 n} } = b_n\Big(\cos \tet, {|p| \over |q|}\Big), $$
}
\noindent
having put
\beq b_n : [-1,1] \times [0,+\infty) \vain [0,+\infty)~, \label{debn} \feq
$$ (c, \tu) \mapsto b_n(c,\tu) := \dd{(1-c^2)(1 + 2 c \tu + \tu^2)^{n}  \over 1 + \tu^{2 n}}~. $$
One checks by elementary tools that
\beq \sup_{c \in [-1,1], \, \tu \in [0,+\infty)} b_n(c,\tu) = b_n\Big({n \over n+2}, 1\Big) =
~{2^{2 n+1} (n+1)^{n+1} \over (n+2)^{n+2}}. \feq
Returning to Eq. \rref{eqgivea}, and writing
$b_n\Big(\cos \tet, \dd{|p| \over |q|}\Big) \leqs
\dd{2^{2 n+1} (n+1)^{n+1} \over (n+2)^{n+2}}$ we obtain the thesis \rref{onthe}. \fine
\begin{prop}
\label{lemz}
\textbf{Lemma.} Let $\nu \in (d,+\infty)$. For any $\ro \in (2 \sqrt{d},+\infty)$, one has
\beq \sum_{h \in \Zd, |h| \geqs \ro} {1 \over |h|^{\nu}}~\leqs
{2 \pi^{d/2} \over \Gamma(d/2)}
\sum_{i=0}^{d-1} \left( \barray{c} d - 1 \\ i \farray \right) {d^{d/2-1/2-i/2} \over (\nu - 1 - i)
(\ro - 2 \sqrt{d})^{\nu - 1 - i}}~. \label{desnu}\feq
\end{prop}
\textbf{Proof.} This is just Lemma C.2 of  \cite{tre} (with the variable $\lan$ of the cited reference
related to $\ro$ by $\lan = \ro - 2 \sqrt{d}$). \fine
\begin{prop}
\textbf{Lemma.}
Let $\ro \in (1,+\infty)$ and $\varphi: [1,\ro) \vain \reali$.
Then, for each $k \in \reali^d$,
\beq
\sum_{h \in \Zd_0, |h| < \ro} (h \sc k)^2  \varphi(|h|) =
{|k|^2 \over d} \sum_{h \in \Zd_0, |h| < \ro} |h|^2 \varphi(|h|)~.
\label{idew} \feq
\end{prop}
\textbf{Proof.} We reexpress the left
hand side of \rref{idew} writing
$(h \sc k)^2 = (\sum_{r=1}^d h_r k_r)$ $\times (\sum_{s=1}^d h_s k_s) = \sum_{r,s=1}^d h_r h_s k_r k_s$,
which gives
\beq \sum_{h \in \Zd_0, |h| < \ro} (h \sc k)^2  \varphi(|h|)
= \sum_{r,s=1}^d \YY_{r s} k_r k_s
\qquad \YY_{r s} :=  \sum_{h \in \Zd_0, |h| < \ro} h_r h_s \varphi(|h|)~.
\label{espw} \feq
From the definition of $\YY_{r s}$, one easily checks that
\beq \YY_{r s} = 0~\mbox{for $r \neq s$},
~~\YY_{1 1} = \YY_{2 2} = ... = \YY_{d d}~. \feq
By the second of the above statements, for each $r \in \{1,...,d\}$ we have
$$ \YY_{r r} = {1 \over d} \sum_{s=1}^d \YY_{s s} =
{1 \over d} \!\!\sum_{h \in \Zd_0, |h| < \ro} \sum_{s=1}^d h^2_s \varphi(|h|)
= {1 \over d} \!\! \sum_{h \in \Zd_0, |h| < \ro} |h|^2 \varphi(|h|)~; $$
in conclusion,
\beq \YY_{r s} = {\delta_{r s} \over d} \!\!
\sum_{h \in \Zd_0, |h| < \ro} |h|^2 \varphi(|h|)~ \qquad (r,s=1,...,d)~. \feq
Inserting this result into the first equality \rref{espw}, we obtain
the thesis \rref{idew}. \fine
\begin{prop}
\label{deden}
\textbf{Definition.} Let us introduce the domain
\beq \Do := \{ (c,\xi) \in \reali^2~|~c \in [-1,1], \, \xi \in [0,+\infty), \, (c,\xi) \neq (1,1) \}~;
\feq
furthermore, let $n \in \reali$. \parn
(i) $E_n$ is the $C^\infty$ function defined as follows:
\beq E_n : \Do \vain [0,+\infty)~,
\quad (c, \ep) \mapsto E_n(c,\ep) := {1 - c^2 \over (1 - 2 c \ep + \ep^2)^{n+1}}~. \label{den} \feq
(ii) For $\ell = 0,1,2,...$, we put
\beq E_{n \ell} :  [-1,1] \vain \reali~, \qquad c \mapsto E_{n \ell}(c) :=
{1 \over \ell!} {\partial^\ell E_n \over \partial \ep^\ell}(c,0)~. \label{denel} \feq
(iii) For $t = 1,2,...$,
\beq R_{n t} : \Do \vain \reali~, \feq
is the unique $C^\infty$ function such that, for all $(c, \ep) \in \Do$,
\beq E_n(c, \ep) = \sum_{\ell=0}^{t-1} E_{n \ell}(c) \ep^\ell + R_{n t}(c, \ep) \ep^t~.
\label{eqrt} \feq
(iv) For $t = 1,2,...$, we put
\beq \mu_{n t} := \min_{c \in [-1,1], \, \ep \in [0,1/2]} R_{n t}(c, \ep)~,
\qquad M_{n t} := \max_{c \in [-1,1], \, \ep \in [0,1/2]} R_{n t}(c, \ep)~.
\label{mnt} \feq
\end{prop}
\begin{rema}
\textbf{Remarks.} (i) Of course, $E_n$ could be defined on a domain
larger than $\Do$; this is not relevant for our purposes. \parn
(ii) Some calculations give
\beq E_{n 0}(c) = 1 - c^2~, \qquad E_{n 1}(c) = (2 n + 2) (c - c^3)~, \label{en02} \feq
$$ E_{n 2}(c)  = -(n+ 1) + (2 n^2 + 7 n + 5) c^2 - (2 n^2 + 6 n + 4) c^4~, $$
etc. \parn
(iii) In general, $E_{n \ell}$ is a polynomial in $c$ of degree $\ell + 2$;
this polynomial is \textsl{even} ($E_{n \ell}(-c) = E_{n \ell}(c)$) for even $\ell$,
and \textsl{odd} ($E_{n \ell}(-c) = - E_{n \ell}(c)$) for odd $\ell$. \parn
(iv) Eq. \rref{eqrt} indicates that $R_{n t}(c, \ep) \ep^t$ is the
reminder in the Taylor expansion of $E_{n}(c, \ep)$ at order $t$ in $\ep$, about
the point $\ep=0$. For the practical
computation of $R_{n t}$ one can note that Eq. \rref{eqrt} (and Taylor's formula) imply
\beq R_{n t}(c, \ep) = \left\{ \barray{ll}
\ep^{-t}  \left(E_n(c, \ep) - \sum_{\ell=0}^{t-1} E_{n \ell}(c) \ep^\ell\right)
& \mbox{if $\ep \neq 0$,} \\
E_{n t}(c) & \mbox{if $\ep=0$}. \farray \right. \label{repnt} \feq
(v) The minimum $\mu_{n t}$ and the maximum $M_{n t}$ in \rref{mnt} exist,
since we consider the continuous function $R_{n t}$ on a compact domain.
For specific values of $n$ and $t$, these can be evaluated numerically
starting from the representation \rref{repnt} of $R_{n t}$. In this way we
obtain, for example, the values
\parn
\vbox{
\beq \mu_{2 6} = -22.720...\,,~ M_{2 6} = 73.835...\,;~~
\mu_{3 6} = -61.239...\,,~M_{3 6} =  410.74...\,; \label{m26} \feq
$$~\mu_{4 6} = -135.89...\,,~
M_{4 6} = 1832.8...\,;~~\mu_{5 6} = -264.44...\,,~M_{5 6} = 7252.9...\,;~$$
$$ \mu_{10, 6} = -2582.5...\,,~ M_{10, 6} = 4.6371... \times 10^6 \,,
$$
}
\noindent
recorded here for subsequent use. \parn
(vi) For arbitrary $n$ and $t$, Eqs. \rref{eqrt} \rref{mnt} imply
\beq \sum_{\ell=0}^{t-1} E_{n \ell}(c) \ep^\ell + \mu_{n t} \ep^t \leqs
E_n(c,\ep) \leqs \sum_{\ell=0}^{t-1} E_{n \ell}(c) \ep^\ell + M_{n t} \ep^t~\label{enfor} \feq
$$ \mbox{for $(c,\ep) \in [-1,1] \times [0,1/2]$}. \qquad \qquad \square $$
\end{rema}
Hereafter we present
a lemma about the function
$k, h \mapsto  \dd{|k|^{2 n} |h \we k|^2 \over |h|^{2 n + 2} |k-h|^{2 n + 2}}$,
appearing in the definition \rref{kknd} of $\KK_{n}$; as indicated by the Lemma,
this is related to the function $E_n$ of Definition \ref{deden} and to its Taylor expansion.
\begin{prop}
\textbf{Lemma.} Let $h, k \in \reali^d \setminus \{0\}$, $h \neq k$, and
let $\tet(h,k) \equiv \tet$ be the convex angle between them. Furthermore,  let $n \in \reali$;
then the following holds. \parn
(i) One has
\beq
{|k|^{2 n} |h \we k|^2 \over |h|^{2 n + 2} |k-h|^{2 n + 2}} =
{1 \over |h|^{2 n}} E_n \Big(\cos \tet, {|h| \over |k|} \Big)~. \label{funf} \feq
(ii) Let $|k| \geqs 2 |h|$. For $t \in \{1,2,...\}$, Eq. \rref{funf} implies
\beq
\sum_{\ell=0}^{t-1} {E_{n \ell}(\cos \tet) \over |h|^{2 n - \ell} |k|^\ell} +
{\mu_{n t} \over  |h|^{2 n - t} |k|^t}  \leqs {|k|^{2 n} |h \we k|^2 \over |h|^{2 n + 2} |k-h|^{2 n + 2}}
 \label{eqtay} \feq
$$
\leqs
\sum_{\ell=0}^{t-1} {E_{n \ell}(\cos \tet) \over |h|^{2 n - \ell} |k|^\ell} +
{M_{n t} \over  |h|^{2 n - t} |k|^t} $$
(note that $\cos \tet = \vh \sc \vk$).
\end{prop}
\textbf{Proof.} (i) Writing $|h \we k|^2 = |h|^2 |k|^2 (1 - \cos^2 \tet)$
and $|k-h|^2$ $= |k|^2 - 2 |k| |h| \cos \tet$ $+ |h|^2$, we readily obtain
$$ {|k|^{2 n} |h \we k|^2 \over |h|^{2 n + 2} |k-h|^{2 n + 2}}
= {1 \over |h|^{2 n}} \, {1 - \cos^2 \tet \over (1 - 2 \cos \tet |h|/|k| + |h|^2/|k|^2)^{n + 1}}~;
$$
comparing this with the definition \rref{den} of $E_n$, we obtain the thesis \rref{funf}. \parn
(ii) It suffices to use Eq. \rref{funf} and the inequalities \rref{enfor}, with
$c := \cos \tet$ and $\ep := |h|/|k|$;
note that $0 \leqs \ep \leqs 1/2$ due to the assumption $|k| \geqs 2 |h|$. \fine
To conclude, we introduce some variants of the polynomials
$E_{n \ell}$, to be used in the sequel.
\begin{prop}
\label{dendel}
\textbf{Definition.} For $\ell = 0,2,...,$, $\Ed_{n \ell d} \equiv \Ed_{n \ell}$ are the polynomials $E_{n \ell}$
of Eq. \rref{denel}, where the term $c^2$ has been replaced with $1/d$.
\end{prop}
\begin{rema}
\textbf{Example.} The expressions of $E_{n 0}$, $E_{n 2}$ in \rref{en02} imply
\beq \Ed_{n 0}(c) = \mbox{const.} = 1 - {1 \over d}~, \label{en02d} \feq
$$ \Ed_{n 2}(c)  = -(n+ 1) + {2 n^2 + 7 n + 5 \over d}  - (2 n^2 + 6 n + 4) c^4~. $$
\end{rema}
\noindent
\section{The function $\boma{\KK_{n}}$}
\label{appek}
Throughout the Appendix, $n \in ({d/2}, + \infty)$.
For $k \in \Zd_0$, let us recall the definition \rref{kknd}
$$ \KK_{n}(k) := |k|^{2 n} \sum_{h \in \Zd_{0 k}} {|h \we k|^2 \over |h|^{2 n + 2} |k-h|^{2 n + 2}}
\in (0,+\infty)~. $$
\begin{prop}
\label{prokknd}
\textbf{Proposition.}
Let us fix a  ''cutoff''
\beq \ro \in (2 \sqrt{d},+\infty)~; \label{cuto} \feq
then, the following holds (with the functions and quantities
$\Km_{n}$, $\dK_{n}$,... mentioned in the sequel
depending parametrically on $d$ and $\ro$: $\Km_{n}(k) \equiv \Km_{n d}(\ro, k)$,
$\dK_{n} \equiv \dK_{n d}(\ro)$,...). \parn
(i) The function $\KK_{n}$
can be evaluated using the inequalities
\beq \Km_{n}(k) < \KK_{n}(k) \leqs \Km_{n}(k) + \dK_{n}
~~\mbox{for all $k \in \Zd_0$}~. \label{dkmnd} \feq
Here
\beq \Km_{n}(k) := |k|^{2 n} \sum_{h \in \Zd_{0 k}, |h| < \ro \op |k-h| < \ro}
 {|h \we k|^2 \over |h|^{2 n + 2} |k-h|^{2 n + 2}}~; \label{deakm} \feq
this can be reexpressed as
\beq \Km_{n}(k) = |k|^{2 n} \sum_{h \in \Zd_{0 k}, |h| < \ro}
\big[1 + \teta(|k-h| - \ro)\big] \, {|h \we k|^2  \over |h|^{2 n + 2} |k-h|^{2 n + 2}}
\label{kmnd} \feq
(with $\teta$ as in Definition \ref{versk}).
If $|k| \geqs 2 \ro$, in Eq. \rref{kmnd} one can replace $\Zd_{0 k}$ with $\Zd_0$
and $\teta(|k-h| - \ro)$ with $1$. \parn
Finally,
\beq \dK_{n} := {2^{2 n + 3} \pi^{d/2} (n+1)^{n+1} \over \Gamma(d/2) (n+2)^{n+2}}
\sum_{i=0}^{d-1} \binom{d-1}{i} {d^{d/2 - 1/2 - i/2} \over (2 n - i - 1)
(\ro - 2 \sqrt{d})^{2 n - i - 1}}~. \label{dedek} \feq
(ii) As in Eq. \rref{refper}, consider the reflection operators
$R_r$
($r=1,...,d$)
and the permutation operators $P_{\si}$ ($\si$ a permutation of $\{1,...,d\}$). Then
\beq \Km_{n}(R_r k) = \Km_{n}(k)~, \quad \Km_{n}(P_\sigma k) = \Km_{n}(k)~
\qquad \mbox{for each $k \in \Zd_0$} \label{claima} \feq
(so, the computation of $\Km_{n}(k)$
can be reduced to the case $k_1 \geqs k_2 \geqs ... \geqs k_d \geqs 0$). \parn
(iii) Let $t \in \{2,4,6,...\}$. One has
$$ \Zz_{n} +
\sum_{\ell=2,4,....,t-2} {\QQ_{n \ell}(\vk) \over |k|^{\ell}} + {v_{n t} \over |k|^t} \leqs
\Km_{n}(k) \leqs
\Zz_{n} +
\sum_{\ell=2,4,....,t-2} {\QQ_{n \ell}(\vk) \over |k|^{\ell}} + {V_{n t} \over |k|^t}
$$
\beq \mbox{for $k \in \Zd_0$, $|k| \geqs 2 \ro$}. \label{takm} \feq
In the above $\sum_{\ell=2,...,t-2} ... := 0$ if $t=2$, and $\vk \in \Sd$
is the versor of $k$ (see Definition \ref{versk}); furthermore,
\beq \Zz_{n} := 2 \left(1 - {1 \over d}\right)
\sum_{h \in \Zd_0, |h| < \ro} {1 \over |h|^{2 n}}~; \label{dekinf} \feq
\beq v_{n t} := 2 \mu_{n t} \hspace{-0.5cm} \sum_{h \in \Zd_0, |h| < \ro} \hspace{-0.4cm} {1 \over |h|^{2 n - t}},
~~V_{n t} := 2 M_{n t} \hspace{-0.5cm} \sum_{h \in \Zd_0, |h| < \ro} \hspace{-0.4cm} {1 \over |h|^{2 n - t}}~~
\mbox{($\mu_{n t}, M_{n t}$ as in \rref{mnt})}~; \label{wnt} \feq
\beq \QQ_{n \ell} : \Sd \vain \reali,~
u \mapsto \QQ_{n \ell}(u) := 2 \hspace{-0.4cm} \sum_{h \in \Zd_0, |h| < \ro} \hspace{-0.4cm}
{\Ed_{n \ell}(u \sc \vh) \over |h|^{2 n - \ell}}~
(\mbox{$\Ed_{n \ell}$ as in Definition \ref{dendel}}). \label{deqnel} \feq
For each $\ell$, $\QQ_{n \ell}$ is a polynomial function on $\Sd$; setting
\beq q_{n \ell} := \min_{u \in \Sd} \QQ_{n \ell}(u)~, \qquad
Q_{n \ell} := \max_{u \in \Sd} \QQ_{n \ell}(u)~, \label{mqnel} \feq
one infers from \rref{takm} that
$$ \Zz_{n} +
\sum_{\ell=2,4,....,t-2} {q_{n \ell} \over |k|^{\ell}} + {v_{n t} \over |k|^t} \leqs
\Km_{n}(k) \leqs \Zz_{n} +
\sum_{\ell=2,4,....,t-2} {Q_{n \ell} \over |k|^{\ell}} + {V_{n t} \over |k|^t}$$
\beq
\qquad \mbox{for $k \in \Zd_0$, $|k| \geqs 2 \ro$}~. \label{tkm} \feq
These facts imply
\beq \Km_{n}(k) \vain \Zz_n \qquad \mbox{for $k \vain \infty$}~. \label{liminf} \feq
(iv) Items (i) and (iii) imply
\beq \sup_{k \in \Zd_0} \Km_{n}(k)  \leqs
\sup_{k \in \Zd_0} \KK_{n}(k) \leqs \Big(\sup_{k \in \Zd_0} \Km_{n}(k) \Big) +
\delta \Km_{n} < + \infty~. \label{supk} \feq
\end{prop}
\textbf{Proof.} We fix a cutoff
$\ro$ as in \rref{cuto}, and
proceed in several steps. More precisely
Steps 1-5 give proofs of statements (i)(ii), while
Steps 6-8 prove statements (iii)(iv).
The assumption \rref{cuto} $\ro > 2 \sqrt{d}$ is essential in Step 3.
\parn
\textsl{Step 1. One has
\beq \KK_{n}(k)= \Km_{n}(k) + \DK_{n}(k)~\quad \mbox{for all $k \in \Zd_0$}~,
\label{decomp} \feq
with
$\Km_{n}(k) := |k|^{2 n} \dd{\sum_{h \in \Zd_{0 k}, |h| < \ro \op |k-h| < \ro}
 {|h \we k|^2 \over |h|^{2 n + 2} |k-h|^{2 n + 2}}}$,
as in \rref{deakm}, and
\beq \DK_{n}(k) := |k|^{2 n} \sum_{h \in \Zd_0, |h| \geqs \ro, |k-h| \geqs \ro}
{|h \we k|^2 \over |h|^{2 n + 2} |k-h|^{2 n + 2}} \in (0,+\infty)~. \label{deakp} \feq}
The above decomposition follows noting that $\Zd_{0 k}$ is the disjoint union of
the domains of the sums defining $\Km_n(k)$ and $\DK_{n}(k)$.
$\Km_{n}(k)$ is finite, involving finitely many summands; $\DK_{n}(k)$ is
finite as well, since we know that $\KK_{n}(k) < +\infty$. \parn
\textsl{Step 2. For each $k \in \Zd_0$, one has the representation \rref{kmnd}
$$ \Km_{n}(k) = |k|^{2 n} \sum_{h \in \Zd_{0 k}, |h| < \ro}
\big[1 + \teta(|k-h| - \ro) \big]~
{|h \we k|^2 \over |h|^{2 n + 2} |k-h|^{2 n + 2}}~. $$
If $|k| \geqs 2 \ro$, in Eq.\rref{kmnd}
one can replace $\Zd_{0 k}$ with  $\Zd_{0}$
and $\teta(|k-h| - \ro)$ with $1$}. \parn
To prove \rref{kmnd} we start from the definition \rref{deakm} of $\Km_n(k)$,
and reexpress the sum therein using
Eq. \rref{tesif}, with $f(h) \equiv f_k(h) :=
\dd{|h \we k|^2 \over |h|^{2 n + 2} |k-h|^{2 n + 2}}$. Noting that
$f(k-h) = f(h)$, we can write
$$ f(h) + \teta(|k-h| - \ro) f(k-h) = \big[1 + \teta(|k-h| - \ro)\big] {|h \we k|^2 \over |h|^{2 n + 2}
|k-h|^{2 n + 2}}~, $$
and obtain Eq. \rref{kmnd}. Now, assume
$|k| \geqs 2 \ro$; then, for all $h \in \Zd_0$ with $|h| < \ro$ one has
\beq |k - h| \geqs |k| - |h| > \ro~; \feq
this implies $h \neq k$ (i.e., $h \in \Zd_{0 k}$) and $\theta(|k-h| - \ro) = 1$, two facts which
justify the replacements indicated above in \rref{kmnd}.
\parn
\textsl{Step 3. For each $k \in \Zd_0$, one has
\beq 0 <  \DK_{n}(k) \leqs \dK_{n}~, \label{boundk} \feq
with $\dK_n$ as in Eq. \rref{dedek}.}
The obvious relation $0 < \DK_{n}(k)$ was already noted; in the sequel we prove that
$\DK_{n}(k) \leqs \dK_{n}$.
To show this, we note the following: for each $h$ in the sum \rref{deakp},
one can write
\parn
\vbox{
\beq  |h \we k|^2 |k|^{2 n} = |h \we (k-h)|^2 |h + (k-h)|^{2 n}  \label{rem2} \feq
$$ \leqs B_n (|h|^{2 n + 2} |k-h|^2 + |h|^2 |k-h|^{2 n + 2} )~, $$ }
\noindent
where the last passage depends on the inequality \rref{onthe}, applied
with $p=h$ and $q=k-h$; for the sake of brevity, we have put
\beq B_n := {2^{2 n+1} (n+1)^{n+1} \over (n+2)^{n+2}}~\label{defbn}~. \feq
Eqs. \rref{deakp} \rref{rem2} give
\beq \DK_{n}(k) \leqs
B_n \, \Big( \sum_{h \in \Zd_0, |h| \geqs \ro, |k-h| \geqs \ro} {1 \over |k-h|^{2 n}} +
\sum_{h \in \Zd_0, |h| \geqs \ro, |k-h| \geqs \ro} {1 \over |h|^{2 n}} \Big)~. \feq
The domain of the above two sums is contained in each one of the sets
$\{ h \in \Zd_0~|$ $~|k-h | \geqs \ro \}$ and
$\{ h \in \Zd_0~|~|h | \geqs \ro \}$; so,
$$ \DK_{n}(k) \leqs B_n \, \Big( \sum_{h \in \Zd_0, |k-h | \geqs \ro}
{1 \over |k-h|^{2 n}} + \sum_{h \in \Zd_0, |h| \geqs \ro} {1 \over |h|^{2 n}} \Big)~. $$
Now, the change of variable $h \mapsto k - h$ in the first sum shows that it is equal to
the second one, so
\beq \DK_{n}(k) \leqs 2 B_n \sum_{h \in \Zd_0, |h | \geqs \ro}
{1 \over |h|^{2 n}}~. \label{eef} \feq
Finally, Eq. \rref{eef} and Eq. \rref{desnu} with $\nu = 2 n$ give
$$ \DK_{n}(k) \leqs  ~
{4 \pi^{d/2} B_n \over \Gamma(d/2)}
\sum_{i=0}^{d-1} \left( \barray{c} d - 1 \\ i \farray \right) {d^{d/2-1/2-i/2} \over (2 n - i - 1)
(\ro - 2 \sqrt{d})^{2 n - i - 1}}~; $$
the right hand side of this inequality is $\dK_n$ defined by \rref{dedek}, as seen immediately
using the definition \rref{defbn} of $B_n$.
\parn
\textsl{Step 4. One has the equalities \rref{claima}
$\Km_{n}(R_r k) = \Km_{n}(k)$, $\Km_{n}(P_\sigma k) = \Km_{n}(k)$,
involving the reflection and permutation operators  $R_r, P_\si$.} The proof
starts from the definition \rref{deakm} of $\Km_n$, and
is very similar to the one employed for the analogous properties
of $\KK_{n}$ (see Eq. \rref{claim} and the subsequent comments). \parn
\textsl{Step 5. One has the inequalities \rref{dkmnd}
$\Km_{n}(k) < \KK_{n}(k) \leqs \Km_{n}(k) + \dK_{n}$.}
These relations follow immediately from the decomposition
\rref{decomp}
$\KK_{n}(k)= \Km_{n}(k) + \DK_{n}(k)$
and from the bounds \rref{boundk} on $\DK_{n}(k)$. \parn
\textsl{Step 6. Let $t \in \{2,4,6,...\}$; one has the inequalities \rref{takm}
for $\Km_n(k)$.}
As an example, we prove the upper bound \rref{takm}
$$ \Km_{n}(k) \leqs \Zz_{n} +
\sum_{\ell=2,4,....,t-2} {\QQ_{n \ell}(\vk) \over |k|^{\ell}} + {V_{n t} \over |k|^t}
\qquad \mbox{for $k \in \Zd_0$, $|k| \geqs 2 \ro$}~. $$
For $k$ as above we can express $\Km_n(k)$
via Eq. \rref{kmnd}, replacing therein $\Zd_{0 k}$ with $\Zd_0$ and
$\theta(|k-h| - \ro)$ with $1$ (see the final statement in Step 2). So,
\beq \Km_{n}(k) = 2 \sum_{h \in \Zd_0, |h | < \ro } \hspace{-0.2cm}
{|k|^{2 n} |h \we k|^2  \over |h|^{2 n + 2} |k-h|^{2 n + 2}}~. \feq
In this expression we insert the upper bound of Eq. \rref{eqtay},
writing therein $\cos \tet = \vh \sc \vk$ (note that
\rref{eqtay} can be used, since $|h|/|k| < \ro/(2 \ro) < 1/2$ for each $h$ in the sum). In this way we obtain
$$ \Km_{n}(k) \leqs
2 \sum_{h \in \Zd_0,
|h | < \ro } \left[\sum_{\ell=0}^{t-1} {E_{n \ell}( \vh \sc \vk ) \over |h|^{2 n - \ell} |k|^\ell} +
{M_{n t} \over  |h|^{2 n - t} |k|^t}~\right] $$
$$ = 2 \sum_{\ell=0}^{t-1}  {1 \over |k|^\ell} \sum_{h \in \Zd_0, |h | < \ro }
{E_{n \ell}( \vh \sc \vk ) \over |h|^{2 n - \ell} } +
{2 M_{n t} \over   |k|^t}  \sum_{h \in \Zd_0, |h | < \ro } {1 \over |h|^{2 n - t}}~; $$
comparing with the definition \rref{wnt}, we see that the last term above
is just $V_{n t}/|k|^t$. Our computation can be summarized in the equation
\beq \Km_{n}(k) \leqs \sum_{\ell=0}^{t-1}  {\QQ_{n \ell}(\vk) \over |k|^\ell} +
{V_{n t} \over   |k|^t}~, \feq
where we have provisionally put
\beq \QQ_{n \ell} : \Sd \vain \reali~, \qquad u \mapsto \QQ_{n \ell}(u) :=
2 \sum_{h \in \Zd_0, |h | < \ro }
{E_{n \ell}(\vh \sc u) \over |h|^{2 n - \ell} }~. \label{deqqnel} \feq
Now, the thesis follows if we prove the following relations:
\beq \QQ_{n 0}(u) = \Zz_n ~\mbox{as in \rref{dekinf}, for all $u \in \Sd$}~; \label{b34} \feq
\beq \QQ_{n \ell}(u) = 0 ~\mbox{for $\ell \in \{1,3,...,t-1\}$ and all $u \in \Sd$}~; \label{b35} \feq
\beq \QQ_{n \ell}(u)~\mbox{is as in \rref{deqnel}, for $\ell \in \{2,4,...,t-2\}$ and all $u \in \Sd$}~.
\label{b36}\feq
To prove \rref{b34}, we proceed as follows,
for any $u \in \Sd$:
recalling that $E_{n 0}(c) = 1 - c^2$ and writing
$\vh = h/|h|$ we get
$$ \QQ_{n 0}(u) =
2 \sum_{h \in \Zd_0, |h | < \ro } {1 \over |h|^{2 n}}
- 2 \sum_{h \in \Zd_0, |h | < \ro } {(h \sc u)^2 \over |h|^{2 n + 2}}~;
$$
on the other hand, the identity \rref{idew} with $k$ replaced by $u$
and $\varphi(|h|) = 1/|h|^{2 n+2}$ gives
$\sum_{h \in \Zd_0, |h | < \ro } {(h \sc u)^2 / |h|^{2 n + 2}} =
(1/d) \sum_{h \in \Zd_0, |h | < \ro } {1/|h|^{2 n}}$. So,
$$ \QQ_{n 0}(u) = 2 \left(1 - {1 \over d} \right) \sum_{h \in \Zd_0, |h | < \ro } {1 \over |h|^{2 n}} =
\Zz_n~, $$
and \rref{b34} is proved. \parn
Let us pass to \rref{b35}. This relation is proved recalling
that, for $\ell$ odd, the function $c \mapsto E_{n \ell}(c)$ is odd as well; this implies
that the general term of the sum \rref{deqqnel} changes its sign under
a transformation $h \mapsto -h$. \parn
Finally, let us prove \rref{b36} for any even $\ell$. In this case we have an
even polynomial
\beq E_{n \ell}(c) = \sum_{j=0,2,...,\ell+2} E_{n \ell j} c^j~, \feq
so \rref{deqqnel} implies
\beq \QQ_{n \ell}(u) = 2 \sum_{j=0,2,...,\ell+2} E_{n \ell j} \sum_{h \in \Zd_0, |h | < \ro }
{(\vh \sc u)^j \over |h|^{2 n - \ell} }~; \label{b40} \feq
in particular, for the $j=2$ term above we have
\beq \sum_{h \in \Zd_0, |h | < \ro } {(\vh \sc u)^2 \over |h|^{2 n - \ell} }~=
\sum_{h \in \Zd_0, |h | < \ro } {(h \sc u)^2 \over |h|^{2 n - \ell + 2} }~=
{1 \over d} \sum_{h \in \Zd_0, |h | < \ro }
{1\over |h|^{2 n - \ell} }~, \label{b41} \feq
where the last passage follows from the identity \rref{idew} (with $k$ replaced by $u$
and $\varphi(|h|) = 1/|h|^{2 n - \ell + 2}$). Eqs. \rref{b40} \rref{b41} imply
\beq
\QQ_{n \ell}(u) = 2 \sum_{j=0,4,6,...,\ell+2} E_{n \ell j} \sum_{h \in \Zd_0, |h | < \ro }
{(\vh \sc u)^j \over |h|^{2 n - \ell} }~ + {2 E_{n \ell 2} \over d}
 \sum_{h \in \Zd_0, |h | < \ro } {1\over |h|^{2 n - \ell} }~.
\label{dacomp} \feq
On the other hand, the Definition \ref{dendel} of $\Ed_{n \ell}$ prescribes
\beq \Ed_{n \ell}(c) = \sum_{j=0,4,6,...,\ell+2} E_{n \ell j} c^j + {E_{n \ell 2} \over d}~; \feq
comparing this with \rref{dacomp}, we conclude
$$
\QQ_{n \ell}(u) = 2 \sum_{h \in \Zd_0, |h | < \ro }
{\Ed_{n \ell}(\vh \sc u) \over |h|^{2 n - \ell} }~~~\mbox{as in \rref{deqnel}}~,
$$
and \rref{b36} is proved.\parn
\textsl{Step 7. Let $t \in \{2,4,6,...\}$. For $\ell \in \{2,4,...,t-2\}$,
the $\QQ_{n \ell}$ are polynomial functions on $\Sd$; considering
their minima $q_{n \ell}$ and maxima $Q_{n \ell}$, one infers from \rref{takm} the inequalities
\rref{tkm}
$$ \Zz_{n} +
\sum_{\ell=2,4,....,t-2} {q_{n \ell} \over |k|^{\ell}} + {v_{n t} \over |k|^t} \leqs
\Km_{n}(k) \leqs \Zz_{n} +
\sum_{\ell=2,4,....,t-2} {Q_{n \ell} \over |k|^{\ell}} + {V_{n t} \over |k|^t}
\quad \mbox{for $|k| \geqs 2 \ro$}~, $$
which imply the relation \rref{liminf}
$\Km_{n}(k) \vain \Zz_n$ for $k \vain \infty$}.
\noindent
The polynomial nature of each function $\QQ_{n \ell}$ follows
from its definition \rref{deqnel} in terms of the polynomial $\Ed_{n \ell}$.
The inequalities \rref{tkm} for $\Km_n(k)$
are obvious; the statement \rref{liminf} follows noting that, in Eq. \rref{tkm},
both the lower and the upper bound for $\Km_n(k)$ tend to $\Zz_n$ for $k \vain \infty$. \parn
\textsl{Step 8. Proof of the inequalities \rref{supk}}
$$ \sup_{k \in \Zd_0} \Km_{n}(k) \leqs
\sup_{k \in \Zd_0} \KK_{n}(k) \leqs \Big(\sup_{k \in \Zd_0} \Km_{n}(k) \Big) +
\delta \Km_{n} < + \infty~. $$
The first two of the above inequalities are an obvious consequence of the relations
\rref{dkmnd}
$\Km_{n}(k) < \KK_{n}(k) \leqs \Km_{n}(k) + \dK_{n}$; the third
inequality holds if we show that
\beq \sup_{k \in \Zd_0} \Km_{n}(k) < + \infty~, \feq
and this follows from the existence of a finite $k \vain \infty$ limit for $\Km_n(k)$
(see Step 7).
\fine
\noindent
\section{Appendix. The upper bounds $\boma{K^{+}_{n}}$, for
$\boma{d=3}$  and $\boma{n=2,3,4,5,10}$}
\label{appe345}
Eq. \rref{deup} defines $K^{+}_n$ in terms
of $\sup_{k \in \Zt_0} \KK_n(k)$, or of any upper approximant for this sup.
In all the cases analyzed hereafter, we produce both an upper and a lower
approximant; the lower one is given only to indicate the
uncertainty in our evaluation of $\sup \KK_n$. \salto
\textbf{Some details on the
evaluation of $\boma{\KK_{2}}$ and of its sup.}
Among the examples
presented here, the case of $\KK_{2}$ is
the one requiring more expensive computations. \parn
To evaluate $\KK_{2}$, we apply Proposition \ref{prokknd} with a fairly large cutoff
\beq \ro = 20~;  \feq
thus, we must often sum over the set
$\{ h \in \Zt_0~|~|h| < 20 \}$.
Eq. \rref{dedek} gives
\beq \dK_{2} = 5.6856...~, \label{dk23} \feq
and it remains  to evaluate the function $\Km_{2}$, using directly
the definition \rref{kmnd} or the bounds in Proposition \ref{prokknd}.
\parn
To evaluate $\Km_{2}(k)$, we start from the $k's$ in $\Zt_0$ with $|k| < 2 \ro = 40$.
We use directly the definition \rref{kmnd}
for all such $k$'s ({\footnote{In fact, due to the symmetry properties \rref{claima},
computation of $\Km_2(k)$ can be limited to points $k$
such that $k_1 \geqs k_2 \geqs k_3 \geqs 0$.}}); in this way, we obtain
\beq \max_{k \in \Zt_0, |k| < 40} \Km_{2}(k) = \Km_{2}(9,9,9) = 22.022... \label{c3} \feq
(another result is that $\Km_{2}(k)$ has a small oscillation for $|k|$ between
$10$ and $40$, since $\Km_{2}(k) \geqs 21.563$ for
$k \in \Zt$, $10 < |k| < 40$). \parn
Let us pass to the case $|k| \geqs 40$. Here, our main tool is the upper bound \rref{tkm} with $t=6$;
after some computation, this yields the result
({\footnote{Let us give some supplementary information on the computations
yielding \rref{c5}. The $t=6$ upper bound in \rref{tkm} reads
$$ \Km_2(k) \leqs \Zz_2 + {Q_{2 2} \over |k|^2} + {Q_{2 4} \over |k|^4} + {V_{2 6} \over |k|^6}~, $$
and we must determine the constants $\Zz_2$, etc., appearing therein.
$\Zz_2$ and $V_{2 6}$ are computed directly from the definitions \rref{dekinf} \rref{wnt}
(the second one requiring previous knowledge of $M_{2 6} = 73.835...$, see Eq. \rref{m26}); in
this way one gets $\Zz_2 = 21.204...$ and $V_{2 6} = 1.1794... \times 10^9$.
$Q_{2 2}$ is the maximum of the polynomial function $\QQ_{2 2} : \St \vain \reali$, and
Eq. \rref{deqnel} gives for this function the explicit expression
$$ \QQ_{2 2}(u) = 2904.7...
- 4569.7... \,(u_1^2 u_2^2 + u_1^2 u_3^2 + u_2^2 u_3^2)
- 2349.8... \,(u_1^4 + u_2^4 + u_3^4)~, $$
for all $u \in \St$; one finds
$Q_{2 2} = \QQ_{2 2}(1/\sqrt{3}, 1/\sqrt{3}, 1/\sqrt{3}) = 598.27...\,$. \parn
$Q_{2 4}$ is the maximum of the polynomial function $\QQ_{2 4}: \St \vain \reali$; after
computing this function via Eq. \rref{deqnel}, one gets $Q_{2 4} = 1.1506... \times 10^5$. \parn
After rounding up from above all the numerical outputs, the computation we have just outlined gives
the first inequality \rref{c5} $\Km_2(k) \leqs 21.205 + 598.28 |k|^{-2} +$ etc.,
holding for $|k| \geqs 40$; on the other hand,
$21.205 + 598.28 |k|^{-2} +$ etc. $\leqs 21.912$
for all such $k$'s, which explains the second inequality \rref{c5}}})
\beq \Km_2(k) \leqs 21.205 + {598.28 \over |k|^2} + {1.1507 \times 10^5 \over |k|^4}
+ {1.1795 \times 10^9 \over |k|^6}~ \leqs 21.912 \,. \label{c5} \feq
(For completeness, we also mention that the $t=6$ lower bound in \rref{tkm} and
Eq. \rref{liminf} imply $\inf_{k \in \Zd_0, |k| \geqs 40} \Km_2(k) =
\lim_{k \vain \infty} \Km_2(k) = 21.204...$
({\footnote{First of all, Eq. \rref{liminf} gives
$$ \lim_{k \vain \infty} \Km_2(k) = \Zz_2 $$
where $\Zz_2$ is defined by \rref{dekinf}; from the previous footnote, we
know that $\Zz_2 =  21.204...$. Now, let us pass to the lower bound \rref{tkm},
with $t=6$; after computing all the necessary constants and
rounding up the results from below, we obtain
$$ \Km_2(k) \geqs
\Zz_2 + {554.98  \over |k|^{2}} + {1.1413 \times 10^5 \over |k|^{4}}
- {3.6293 \times 10^8 \over |k|^{6}}~
\quad \mbox{for $k \in \Zt_0$, $|k| \geqs 40$}. $$
On the other hand, one has
$554.98 \,|k|^{-2} + 1.1413 \times 10^5 \,|k|^{-4} ... \geqs 0$
for $|k| \geqs 40$; so, $\inf_{k \in \Zt_0, |k| \geqs 40} \Km_2(k) \geqs \Zz_2$.
It is obvious that $\inf_{k \in \Zt_0, |k| \geqs 40} \Km_2(k) \leqs
\lim_{k \in \Zt_0, k \vain \infty} \Km_2(k)$; the latter equals $\Zz_2$,
thus $\inf_{|k| \geqs 40} \Km_2(k) = \lim_{k \vain \infty} \Km_2(k) = \Zz_2$.}});
by comparison with \rref{c5}, we see that
$\Km_2(k)$ is almost constant for $|k| \geqs 40$.) \parn
The results \rref{c3} \rref{c5} yield
\beq \sup_{k \in \Zt_0} \Km_{2}(k) = \Km_{2}(9,9,9) = 22.022...~. \label{summar} \feq
\parn
\vbox{
\noindent
We now pass to the function $\KK_{2}$; according to \rref{supk} we have
$\sup_{k \in \Zt_0} \Km_{2}(k) \leqs
\sup_{k \in \Zt_0} \KK_{2}(k) \leqs \Big(\sup_{k \in \Zt_0} \Km_{2}(k) \Big) +
\delta \Km_{2}$, and the results \rref{dk23} \rref{summar} give
\beq 22.022 < \sup_{k \in \Zt_0} \KK_{2}(k) < 27.709~. \label{sugiv} \feq
}
\noindent
(The uncertainty on this sup is fairly large, due to the
value of $\dK_2$ in \rref{dk23}; the error $\dK_2$
could be significantly reduced choosing a cutoff $\ro \gg 20$,
but the related computations would be much more expensive.) \parn
\textbf{The upper bound $\boma{K^{+}_{2}}$.}
According to the definition \rref{deup}, we have
\beq K^{+}_{2} = {1 \over (2 \pi)^{3/2}}
\sqrt{\sup_{k \in \Zt_0} \KK_{2}(k)} \quad \mbox{(or any upper approximant for this)}~. \feq
Due to \rref{sugiv}, we can take $K^{+}_2 = (2 \pi)^{-3/2} \sqrt{27.709}\,$;
rounding up to  three digits we can write
\beq K^{+}_{2} = 0.335 \, , \feq
as reported in \rref{bouk}.
\salto
\textbf{Preparing the examples with $\boma{n=3,4,5,10}$.}
To evaluate $\KK_n$ for the cited values of $n$, we apply Proposition \ref{prokknd} with a cutoff
\beq \ro = 10~;  \feq
thus, all sums over $h$ in Proposition \ref{prokknd} are over the set
$\{ h \in \Zt_0~|~|h| < 10 \}$.
\salto
\textbf{Some details on the
evaluation of $\boma{\KK_{3}}$ and of its sup.}
Eq. \rref{dedek} gives
\beq \dK_{3} = 0.45295...~, \label{dk23i} \feq
and it remains  to evaluate the function $\Km_{3}$.
Direct computation of this function from the definition \rref{kmnd},
for all $k$'s of norm $< 20$, gives
\beq \max_{k \in \Zt_0, |k| < 20} \Km_{3}(k) = \Km_{3}(2,1,1) = 25.301...~. \label{cc3} \feq
On the other hand the upper bound in Eq. \rref{tkm},
with $t=6$, gives
\parn
\vbox{
$$ \Km_3(k) \leqs 11.197 + {117.33 \over |k|^2} + {1581.5 \over |k|^4}
+ {3.3994 \times 10^6 \over |k|^6}  \leqs 11.554 $$
\beq \mbox{for $k \in \Zt_0$, $|k| \geqs 20$}~. \qquad \label{cc5} \feq
}
\noindent
(For completeness, we mention that the $t=6$ lower bound in \rref{tkm}
and Eq. \rref{liminf} give
$\inf_{k \in \Zt_0, |k| \geqs 20} \Km_{3}(k) =
\lim_{k \vain \infty} \Km_{3}(k) = 11.196...\,$.
Comparing with \rref{cc5}
we conclude that $\Km_3(k)$ is almost constant for $|k| \geqs 20$.) \parn
Eqs. \rref{cc3} \rref{cc5} yield
\beq \sup_{k \in \Zt_0} \Km_{3}(k) = \Km_{3}(2,1,1) = 25.301...~. \label{summari} \feq
We now pass to the function $\KK_{3}$; according to \rref{supk} we have
$\sup_{k \in \Zt_0} \Km_{3}(k) \leqs
\sup_{k \in \Zt_0} \KK_{3}(k) \leqs \Big(\sup_{k \in \Zt_0} \Km_{3}(k) \Big) +
\delta \Km_{3}$, and the results \rref{dk23i} \rref{summari} give
\beq 25.301 < \sup_{k \in \Zt_0} \KK_{2}(k) < 25.755~. \label{sugivei} \feq
\textbf{The upper bound $\boma{K^{+}_{3}}$.}
According to the definition \rref{deup}, we have
\beq K^{+}_{3} = {1 \over (2 \pi)^{3/2}}
\sqrt{\sup_{k \in \Zt_0} \KK_{3}(k)} \quad \mbox{(or any upper approximant for this)}~. \feq
Due to \rref{sugivei}, we can take $K^{+}_3 = (2 \pi)^{-3/2} \sqrt{25.755}\,$; rounding up to  three digits we can write
\beq K^{+}_{3} = 0.323\,, \feq
as reported in \rref{bouk}. \salto
\textbf{Some details on the
evaluation of $\boma{\KK_{4}}$ and of its sup.}
Eq. \rref{dedek} gives
\beq \dK_{4} = 0.021561...~, \label{dk23ib} \feq
and it remains  to evaluate the function $\Km_{4}$.
Direct computation of this function from the definition \rref{kmnd},
for all $k$'s of norm $< 20$, gives
\beq \max_{k \in \Zt_0, |k| < 20} \Km_{4}(k) = \Km_{4}(2,1,0) = 48.038...~. \label{ccc3} \feq
On the other hand the upper bound in Eq. \rref{tkm},
with $t=6$, gives
\beq \Km_4(k) \leqs 9.2611 + {137.37 \over |k|^2} + {629.55 \over |k|^4}
+ {4.2612 \times 10^5 \over |k|^6} \leqs 9.6152 \label{ccc5} \feq
$$ \mbox{for $k \in \Zt_0$, $|k| \geqs 20$}~.  $$
(For completeness we mention that the $t=6$ lower bound in \rref{tkm}
implies $\Km_4(k) \geqs  9.2380$ for $k \in \Zt_0$, $|k| \geqs 20$
({\footnote{More precisely: the $t=6$ lower bound in \rref{tkm} gives
$\Km_4(k) \geqs 9.2610 - 10.098 \,|k|^{-2} + 446.33 \,|k|^{-4}
- 3.1595 \times 10^4 \,|k|^{-6}$ for $k \in \Zt_0$, $|k| \geqs 20$;
for $k$ in the same range, one has $9.2610 - 10.098 \,|k|^{-2} + ... \geqs
9.2380$.}}), while
Eq. \rref{liminf} gives
$\lim_{k \vain \infty} \Km_4(k) = 9.2610...\,.$) \parn
Eqs. \rref{ccc3} \rref{ccc5} yield
\beq \sup_{k \in \Zt_0} \Km_{4}(k) = \Km_{4}(2,1,0) = 48.038...~. \label{summarib} \feq
We now pass to the function $\KK_{4}$; according to \rref{supk} we have
$\sup_{k \in \Zt_0} \Km_{4}(k) \leqs
\sup_{k \in \Zt_0} \KK_{4}(k) \leqs \Big(\sup_{k \in \Zt_0} \Km_{4}(k) \Big) +
\delta \Km_{4}$, and the results \rref{dk23ib} \rref{summarib} give
\beq 48.038 < \sup_{k \in \Zt_0} \KK_{4}(k) < 48.061~. \label{sugiveib} \feq
\textbf{The upper bound $\boma{K^{+}_{4}}$.}
According to the definition \rref{deup}, we have
\beq K^{+}_{4} = {1 \over (2 \pi)^{3/2}}
\sqrt{\sup_{k \in \Zt_0} \KK_{4}(k)} \quad \mbox{(or any upper approximant for this)}~. \feq
Due to \rref{sugiveib}, we can take $K^{+}_4 = (2 \pi)^{-3/2} \sqrt{48.061}\,$; rounding up to  three digits we can write
\beq K^{+}_{4} = 0.441 \,, \feq
as reported in \rref{bouk}. \salto
\textbf{Some details on the
evaluation of $\boma{\KK_{5}}$ and of its sup.}
Eq. \rref{dedek} gives
\beq \dK_{5} = 0.0012414...~, \label{dk23ic} \feq
and it remains  to evaluate the function $\Km_{5}$.
Direct computation of this function from the definition \rref{kmnd},
for all $k$'s of norm $< 20$, gives
\beq \max_{k \in \Zt_0, |k| < 30} \Km_{5}(k) = \Km_{5}(2,1,0) = 106.99...~ . \label{cccc3} \feq
On the other hand the upper bound in Eq. \rref{tkm},
with $t=6$, gives
$$ \Km_5(k) \leqs 8.5682 + {186.23. \over |k|^2} + {919.89 \over |k|^4}
+ {2.2152 \times 10^5 \over |k|^6} \leqs 9.0430 $$
\beq \mbox{for $k \in \Zt_0$, $|k| \geqs 20$}~. \label{cccc4} \feq
(For completeness we mention that the $t=6$ lower bound in \rref{tkm}
implies
$\Km_5(k) \geqs 8.4974$ for $k \in \Zt_0$, $|k| \geqs 20$, while
Eq. \rref{liminf} gives
$\lim_{k \vain \infty} \Km_5(k) = 8.5681...\,.$) \parn
Eqs. \rref{cccc3} \rref{cccc4} yield
\beq \sup_{k \in \Zt_0} \Km_{5}(k) = \Km_{5}(2,1,0) = 106.99...~. \label{summaric} \feq
We now pass to the function $\KK_{5}$; according to \rref{supk} we have
$\sup_{k \in \Zt_0} \Km_{5}(k) \leqs
\sup_{k \in \Zt_0} \KK_{5}(k) \leqs \Big(\sup_{k \in \Zt_0} \Km_{5}(k) \Big) +
\delta \Km_{5}$, and the results \rref{dk23ic} \rref{summaric} give
\beq 106.99 < \sup_{k \in \Zt_0} \KK_{5}(k) < 107~. \label{sugiveic} \feq
\textbf{The upper bound $\boma{K^{+}_{5}}$.}
According to the definition \rref{deup}, we have
\beq K^{+}_{5} = {1 \over (2 \pi)^{3/2}}
\sqrt{\sup_{k \in \Zt_0} \KK_{5}(k)} \quad \mbox{(or any upper approximant for this)}~. \feq
Due to \rref{sugiveic}, we can take $K^{+}_5 = (2 \pi)^{-3/2} \sqrt{107}\,$; rounding up to  three digits we can write
\beq K^{+}_{5} = 0.657 \,, \feq
as reported in \rref{bouk}. \salto
\textbf{Some details on the
evaluation of $\boma{\KK_{10}}$ and of its sup.}
Eq. \rref{dedek} gives
\beq \dK_{10} = 2.1401... \times 10^{-9}~, \label{dk23id} \feq
and it remains  to evaluate the function $\Km_{10}$.
Direct computation of this function from the definition \rref{kmnd},
for all $k$'s of norm $< 20$, gives
\beq \max_{k \in \Zt_0, |k| < 20} \Km_{10}(k) = \Km_{10}(2,1,0) = 9556.5...\, . \label{ccccc3} \feq
On the other hand the upper bound in Eq. \rref{tkm},
with $t=6$, gives
\vskip 0.2cm
\parn
\vbox{
\beq \Km_{10}(k) \leqs 8.0159 + {617.05 \over |k|^2} + {9693.2 \over |k|^4}
+ {5.6557 \times 10^7 \over |k|^6} \leqs 10.503 \label{ccccc4} \feq
$$ \mbox{for $k \in \Zt_0$, $|k| \geqs 20$}~. $$
}
(For completeness we mention that
the $t=6$ lower bound in \rref{tkm} implies
$ \Km_{10}(k) \geqs 7.8034$ for $k \in \Zt_0$, $|k| \geqs 20$,
while \rref{liminf} gives $\lim_{k \vain \infty} \Km_{10}(k) = 8.0158...\,.$) \parn
Eqs. \rref{ccccc3} \rref{ccccc4} yield
\beq \sup_{k \in \Zt_0} \Km_{10}(k) = \Km_{10}(2,1,0) =  9556.5...~. \label{summarid} \feq
We now pass to the function $\KK_{10}$; according to \rref{supk} we have
$\sup_{k \in \Zt_0} \Km_{10}(k) \leqs
\sup_{k \in \Zt_0} \KK_{10}(k) \leqs \Big(\sup_{k \in \Zt_0} \Km_{10}(k) \Big) +
\delta \Km_{10}$, and the results \rref{dk23id} \rref{summarid} give
({\footnote{In the MATHEMATICA output for $\Km_{10}(2,1,0)$,
$9556.5$ is followed by a digit different from $9$; so, the digits $9556.5$ do
not change when $\delta \Km_{10}$ is added to this output.}})
\beq \sup_{k \in \Zt_0} \KK_{10}(k) = 9556.5...~. \label{sugiveid} \feq
\textbf{The upper bound $\boma{K^{+}_{10}}$.}
According to the definition \rref{deup}, we have
\beq K^{+}_{10} = {1 \over (2 \pi)^{3/2}}
\sqrt{\sup_{k \in \Zt_0} \KK_{10}(k)} \quad \mbox{(or any upper approximant for this)}~. \feq
Using \rref{sugiveid}, and rounding up to  three digits the final result, we can write
\beq K^{+}_{10} = 6.21 \,, \feq
as reported in \rref{bouk}. \parn
\vskip 0.7cm \noindent
\textbf{Acknowledgments.} This work was partly supported by INdAM and by MIUR, PRIN 2008
Research Project "Geometrical methods in the theory of nonlinear waves and applications".
\vfill \eject 
\end{document}